\documentclass{article}
\usepackage{mathrsfs}
\usepackage{amsmath}
\usepackage{amssymb, amsmath}
\usepackage{graphicx}
\usepackage{fancyhdr}
\catcode`\@=11 \@addtoreset{equation}{section}

\catcode`\@=12
\usepackage{color}

\usepackage[colorlinks,
            linkcolor=blue,
            anchorcolor=blue,
           citecolor=red,
           ]{hyperref}
\newtheorem{Theorem}{Theorem}[section]
\newtheorem{Proposition}{Proposition}[section]
\newtheorem{Lemma}{Lemma}[section]
\newtheorem{Corollary}{Corollary}[section]
\newtheorem{Definition}{Definition}[section]

\newtheorem{Remark}{Remark}[section]

\newcommand{\newcom}{\newcommand}
\newcommand{\bTheorem}[1]{
\begin{Theorem} \label{T#1} }
\newcommand{\eT}{\end{Theorem}}

\newcommand{\bProposition}[1]{
\begin{Proposition} \label{P#1}}
\newcommand{\eP}{\end{Proposition}}

\newcommand{\bLemma}[1]{
\begin{Lemma} \label{L#1} }
\newcommand{\eL}{\end{Lemma}}

\newcommand{\bCorollary}[1]{
\begin{Corollary} \label{C#1} }
\newcommand{\eC}{\end{Corollary}}

\newcommand{\beq}{\begin{equation}}
\newcommand{\eeq}{\end{equation}}
\newcom{\ben}{\begin{eqnarray}}
\newcom{\een}{\end{eqnarray}}
\newcom{\beno}{\begin{eqnarray*}}
\newcom{\eeno}{\end{eqnarray*}}
\newcom{\bali}{\begin{aligned}}
\newcom{\eali}{\end{aligned}}

\newcommand{\bFormula}[1]{
\begin{equation} \label{#1}}
\newcommand{\eF}{\end{equation}}

\newcommand{\f}{\frac}
\newcommand{\Om}{\Omega}

\newcommand{\p}{\partial}

\newcommand{\vr}{\varrho}

\newcommand{\vt}{\vartheta}
\newcommand{\vu}{\mathbf{u}}

\newcommand{\vv}{\mathbf{v}}

\newcommand{\vc}[1]{{\boldsymbol #1}}
\newcommand{\Div}{{\rm div}}
\newcommand{\Grad}{\nabla}

\newcommand{\dx}{\,{\rm d} x}

\newcommand{\dt}{\,{\rm d} t }

\newcommand{\dS}{{\rm d} S}

\newcommand{\dxdt}{\dx\dt}

\newcommand{\ep}{\varepsilon}

\font\F=msbm10 scaled 1000
\newcommand{\R}{\mbox{\F R}}

\newcommand\Cbox[2]{%
    \newbox\contentbox%
    \newbox\bkgdbox%
    \setbox\contentbox\hbox to \hsize{%
        \vtop{
            \kern\columnsep
            \hbox to \hsize{%
                \kern\columnsep%
                \advance\hsize by -2\columnsep%
                \setlength{\textwidth}{\hsize}%
                \vbox{
                    \parskip=\baselineskip
                    \parindent=0bp
                    #2
                }%
                \kern\columnsep%
            }%
            \kern\columnsep%
        }%
    }%
    \setbox\bkgdbox\vbox{
        \color{#1}
        \hrule width  \wd\contentbox %
               height \ht\contentbox %
               depth  \dp\contentbox
        \color{black}
    }%
    \wd\bkgdbox=0bp%
    \vbox{\hbox to \hsize{\box\bkgdbox\box\contentbox}}%
    \vskip\baselineskip%
}

\begin{document}


\pagestyle{fancy} \lhead{\color{blue}{Compressible MHD equations with general boundary data}} \rhead{\emph{Y.Li, Y.S.Kwon, Y.Sun}}

\title{\bf Weak solutions to a compressible viscous non-resistive MHD equations with general boundary data }

\author{
Yang Li$^\dag$   \,\,\,\,\,\, \,\,\, Young-Sam Kwon$^\ddag$    \,\,\,\,\,\, \,\,\, Yongzhong Sun$^\S$     \\ \\
 $^\dag$School of Mathematical Sciences, \\ Anhui University, Hefei, 230601, People's Republic of China \\ Email: lynjum@163.com \\ \\
 $^\ddag$Department of Mathematics, \\ Dong-A University, Busan, 604-714, Republic of Korea \\ Email:  ykwon@dau.ac.kr \\ \\
 $^\S$Department of Mathematics, \\ Nanjing University, Nanjing, 210093, People's Republic of China \\ Email: sunyz@nju.edu.cn\\ \\
}

\maketitle

{\centerline {\bf Abstract }}
\vspace{2mm}
{This paper is concerned with a compressible MHD equations describing the evolution of viscous non-resistive fluids in piecewise regular bounded Lipschitz domains. Under the general inflow-outflow boundary conditions, we prove
existence of global-in-time weak solutions with finite energy initial data. The present result extends considerably the previous work by Li and Sun [\emph{J. Differential Equations.}, 267 (2019), pp. 3827-3851], where the homogeneous Dirichlet boundary condition for velocity field is treated. The proof leans on the specific mathematical structure of equations and the recently developed theory of open fluid systems. Furthermore, we establish the weak-strong uniqueness principle, namely a weak solution coincides with the strong solution on the lifespan of the latter provided they emanate from the same initial and boundary data. This basic property is expected to be useful in the study of convergence of numerical solutions.
}

\vspace{2mm}


\vspace{2mm}


\tableofcontents

\section{Introduction}
The time evolution of electrically conductive fluids in the presence of external magnetic field is described by the
magnetohydrodynamic equations (`MHD' in short). It is widely applied in astrophysics, thermonuclear
reactions and industry, among many others. Assuming that the fluids are compressible viscous and isentropic, a simplified but well-accepted three-dimensional MHD equations read as
\begin{equation}\label{in1}
\left\{\begin{aligned}
& \p_t  \vr+ \Div(  \vr \vc{u} )=0,\\
&  \p_t(\vr \vc{u}) +\Div (\vr  \vc{u}\otimes \vc{u} ) +\Grad p(\vr)=\Div \mathbb{S}(\Grad \vc{u} ) +(\Grad \times \mathbf{H})\times \mathbf{H} ,   \\
& \p_t \mathbf{H}= \Grad \times (\vc{u} \times \mathbf{H}) - \nu \Grad \times (\Grad \times \mathbf{H}) , \\
& \Div \mathbf{H}=0 . \\
\end{aligned}\right.
\end{equation}
Here, the unknowns $\vr,\vc{u} \in\R^3,\mathbf{H}\in \R^3$ represent the density of fluids, the velocity field and the magnetic field respectively. $p=p(\vr)$ is the scalar pressure relating to the density. To simplify the presentation, we assume the fluids obey the isentropic law, meaning $p(\vr)=\vr^{\gamma}$ with adiabatic exponent $\gamma>1$. $\mathbb{S}(\Grad \vc{u} )$ is the viscous stress tensor given by $\mathbb{S}(\Grad \vc{u} )=\mu (\Grad \vc{u}+\Grad^{t}\vc{u})+\lambda \Div \vc{u} \mathbb{I}  $ with viscosity coefficients $\mu>0,2\mu+\lambda>0$. $\nu \geq 0$ is the resistivity coefficient acting as the magnetic diffusion. In particular, system (\ref{in1}) is usually termed as compressible non-resistive MHD equations when $\nu=0$; it can be applied to model plasmas when the plasmas are strongly collisional, or the resistivity due to these collisions are extremely small. We refer to \cite{AbZh17} and the references therein for more physical explanations of non-resistve MHD system.

There are many interesting and excellent results in the context of incompressible MHD system. For incompressible viscous resistive MHD equations, the results are classical and largely similar to that of incompressible Navier-Stokes system. In practical models, the resistivity coefficient is very small and it is meaningful to consider non-resistive MHD equations. However, this gives rise to extra mathematical challenges. Bardos et al. \cite{BSS88} initiated the global well-posedness for incompressible ideal MHD equations with strong background magnetic field; see also \cite{CaLe1,HXY1} for more recent treatment. Lin et al. \cite{LXZ} initiated the global well-posedness for incompressible viscous non-resistive MHD equations in 2D. We refer to \cite{AbZh17,PZZu1,TW,XZ} for more results on incompressible viscous non-resistive MHD equations in the three-dimensional whole space, three-dimensional periodic domain and an infinite slab.

Analogous to incompressible fluids, there are richful mathematical results to the compressible viscous resistive MHD equations (\ref{in1}). We refer to Ducomet and Feireisl \cite{DF06}, Hu and Wang \cite{HW2} for the existence of global weak solutions. For global small smooth solutions, we refer to Kawashima \cite{Kaw1} and Li et al. \cite{LiXuZw1} with large oscillations and vacuum, to name only a few. A detailed and complete review in this case is impossible and beyond the scope of this article. Remarkably, the diffusion term in the equation of magnetic field provides a nice compactness and regularity properties.

Let us recall several results on the compressible viscous non-resistive MHD equations, which are closely related to our focus. It turns out that zero resistivity coefficient brings lots of extra difficulty in constructing global-in-time solutions. There are satisfactory results in the simplified 1D case. We refer to Jiang and Zhang \cite{JZW} for global well-posedness with large initial data in 1D bounded domain, see also \cite{LS1D} for the large time behavior and \cite{LY2,LY3} for extensions to the Cauchy problem and planar non-resistive fluids. In the presence of strong background magnetic field, there are small data global well-posedness in multi-dimensional space. Wu and Wu \cite{WW} proved the global small solutions in $\R^2$, while Tan and Wang \cite{TW} obtained the global small solutions in a slab $\R^2 \times (0,1)$, see also \cite{LY2022,Z23} for extensions to the heat-conductive fluids. There are also some recent progress on the existence of global-in-time weak solutions to the compressible viscous non-resistive MHD equations in multi-dimensional space. In a special 2D setting, which will be described in details below, we prove the existence of global weak solutions with finite energy initial data under homogeneous Dirichlet boundary conditions, see \cite{LS19,LS21}. Liu and Zhang \cite{LZ21} extended the result of \cite{LS19} allowing non-monotone pressure law and non-constant viscosity coefficient.
Based on the setting proposed in \cite{LS19}, Dong et al. \cite{DWZ23} recently proved the global small solutions in Besov spaces. We refer to a series of recent nice work by Jiang and Jiang \cite{JJ3,JJ4} on magnetic inhibition theory in non-resistive magnetohydrodynamic fluids.

As a consequence, it remains open for the existence of \emph{global-in-time} weak solutions for the compressible viscous non-resistive MHD equations with \emph{general boundary data}. The general boundary data play significant role in both engineering applications and theoretical analysis, such as the convergence and error estimates of numerical solutions.
This is the main motivation of this paper.

Following our previous work \cite{LS19}, we consider the case that the motion of fluids takes place in the plane $\R^2$, and the magnetic field acts on the fluids only in the vertical direction. More precisely, by choosing
\[
\vc{u}=(u^1,u^2,0)(t,x_1,x_2)=:(\vu,0)(t,x_1,x_2)
\]
and
\[
\mathbf{H}=(0,0,b)(t,x_1,x_2),\,\, \vr=\vr(t,x_1,x_2)
\]
in the equations (\ref{in1}) and setting the resistivity coefficient $\nu$ to be zero, one gets
\begin{equation}\label{in2}
\left\{\begin{aligned}
& \p_t  \vr+ \Div(  \vr \vu )=0,\\
&  \p_t(\vr \vu) +\Div (\vr  \vu \otimes \vu ) +\Grad \left( p(\vr)+\f{1}{2}b^2 \right)=\Div \mathbb{S}(\Grad \vu ) ,      \\
& \p_t b + \Div(b \vu)=0  . \\
\end{aligned}\right.
\end{equation}

To simplify the presentation, we suppose that the motion of fluids takes place in a piecewise regular bounded Lipschitz domain $\Omega \subset \R^2$ and the adiabatic exponent $\gamma>1$\footnote{Throughout this paper, the term piecewise regular domain means that the domain enjoys piecewise $C^2$ regularity. It should be stressed that the main result on global existence of weak solutions still holds for a larger class of domains with the so-called \emph{admissible} inflow-outflow boundary, see \cite{KKNN22}.}.
The system (\ref{in2}) is supplemented with the initial conditions
\begin{equation}\label{in3}
\vr|_{t=0}=\vr_0,\,\, b|_{t=0}=b_0,\,\,\vr \vu|_{t=0}=\vr_0 \vu_0,
\end{equation}
together with the general inflow-outflow boundary conditions
\begin{equation}\label{in4}
\vu|_{\p \Omega}= \vu_{B},\,\, \vr|_{ \Gamma^{ \text{in} } } =\vr_{B},\,\, b |_{ \Gamma^{ \text{in} } } =b_{B},
\end{equation}
where
\[
\Gamma^{ \text{in} } :=\{x\in \p \Omega \, |  \, \vu_{B} \cdot \mathbf{n} <0 \}
\]
with $\mathbf{n}$ being the unit outward normal to $\p \Omega$. Due to the nature of the equations (\ref{in2}), no boundary density and boundary magnetic field are prescribed on
\[
\Gamma^{ \text{out} }:= \{x\in \p \Omega \, |  \, \vu_{B} \cdot \mathbf{n} > 0 \}.
\]


Notice that the compressible MHD equations (\ref{in2}) are reminiscent of the compressible two-fluid models, in the sense that the pressure depends on two variables with each satisfying
the continuity equations. Consequently, we would like to recall some known results for two-fluid models. The compressible two-fluid model considered in Vasseur et al. \cite{VWY} read as
\begin{equation}\label{in5}
\left\{\begin{aligned}
& \p_t  \vr+ \Div(  \vr \vu )=0,\\
& \p_t  n+ \Div(  n \vu )=0,\\
&  \p_t((\vr+n) \vu) +\Div (  (\vr+n)  \vu \otimes \vu ) + \Grad p(\vr,n )=\Div \mathbb{S}(\Grad \vc{u} ).    \\
\end{aligned}\right.
\end{equation}
Here $\vr,n$ are the densities of two fluids and $\vu$ in the common velocity field. The pressure takes the form
\begin{equation}\label{in6}
p(\vr,n)=\vr^{\gamma}+n^{\alpha},\,\, \gamma,\alpha>1.
\end{equation}
In the seminal work, Vasseur et al. \cite{VWY} proved the global existence of weak solutions for (\ref{in5})-(\ref{in6}) with finite energy initial data and suitable constraints on the adiabatic exponents, see also Wen \cite{W21} for the recent progress without the domination condition between the initial densities. Based on a new compactness argument, Bresch et al. \cite{BMZ} proved the global existence of weak solutions for the compressible two-fluid model with algebraic pressure closure. Novotn\'{y} and Pokorn\'{y} \cite{NP} considered a fairly general class of pressure-density laws and proved the existence of global weak solutions. It should be stressed that the mentioned results are proved under the homogeneous Dirichlet boundary conditions for velocity field. Very recently, there are some progress on compressible two-fluid models with general inflow-outflow boundary data. Kracmar et al. \cite{KKNN22} proved global existence of weak solutions for a bifluid model for a mixture of two compressible noninteracting fluids with general boundary data. Jin et al. \cite{JKNN21} obtained the existence and stability of dissipative turbulent solutions to a simple bi-fluid model of compressible fluids.
These later two results partially motivate the present work. However, our compressible MHD equations (\ref{in2}) cannot be regarded as a special case of the compressible two-fluid model due to the specific structure.

\section{The main result}


We are now ready to give the definition of weak solutions as follows.
\begin{Definition}\label{def:1}
Let $T\in (0,\infty)$. A triple $(\vr,\vu,b)$ is said to be a weak solution to (\ref{in2}),  with the initial conditions (\ref{in3}) and the inflow-outflow boundary conditions (\ref{in4}) in $(0,T)\times \Om$ if
\begin{itemize}
\item {Regularity class\footnote{ $ f\in C_{\rm weak}([0,T];L^{p}(\Om)  ) $ means that $f\in L^{\infty}(0,T;L^p(\Om))$ and $\int_{\Om} f(t) \phi \dx \in C([0,T])  $ for any $\phi \in L^{p'}(\Om)$ with $\f{1}{p}+\f{1}{p'}=1$. }
\begin{align*}
    &  \vr(t,x),b(t,x)\geq 0  \text{   for a.e.  }(t,x)\in (0,T)\times \Om,          \\
    &   \vr \in C_{\rm weak}([0,T];L^{\gamma}(\Om)  ) \cap L^{\gamma}( 0,T;L^{\gamma} (\p \Om; | \vu_{B} \cdot \mathbf{n}|\dS )    ),         \\
    &   b \in C_{\rm weak}([0,T];L^{2}(\Om)  ) \cap L^{2}( 0,T;L^{2} (\p \Om; | \vu_{B} \cdot \mathbf{n}|\dS )    ),   \\
    &  \vu-\vu_{B} \in L^2(0,T; W_0^{1,2} (\Om;\R^2)),\,\, \vr \vu \in C_{\rm weak}([0,T];L^{  \f{2\gamma}{\gamma+1} }(\Om)  ),\\
    &  \vr |\vu-\vu_{B}|^2   \in L^{\infty}(0,T;L^1(\Om)) ;
\end{align*}
}
\item { The equation of continuity
\begin{align}
    &  \int_0^{\tau} \int_{\Om} (\vr \p_t \phi+ \vr \vu \cdot \Grad \phi) \dxdt-\int_0^{\tau} \int_{ \Gamma^{ \rm{in} } }   \vr_{B}  \vu_{B}\cdot \mathbf{n} \phi \dS \dt   \nonumber     \\
    &   \quad  = \int_0^{\tau} \int_{ \Gamma^{ \rm{out} } }   \vr  \vu_{B}\cdot \mathbf{n} \phi \dS \dt + \left[  \int_{\Om} \vr \phi \dx  \right]_{t=0}^{t=\tau}    \label{equ-conti}
\end{align}
      for any $\tau \in [0,T]$ and any $\phi\in C^{1}( [0,T] \times \overline{\Om})$;
}
\item{ The balance of momentum
\begin{align}
    &  \int_0^{\tau} \int_{\Om} \left(   \vr \vu \cdot \p_t \vc{\phi} +\vr \vu \otimes \vu :\Grad \vc{\phi}+\left( p(\vr)+\f{1}{2}b^2 \right)\Div \vc{\phi} -\mathbb{S}(\Grad \vu):\Grad \vc{\phi}                     \right) \dxdt    \nonumber     \\
    &   \quad  =  \left[  \int_{\Om} \vr \vu \cdot \vc{\phi} \dx  \right]_{t=0}^{t=\tau}        \label{balan-momentu}
\end{align}
      for any $\tau \in [0,T]$ and any $\vc{\phi} \in C_c^{1}( [0,T] \times \Om;\R^2)$;
}
\item { The equation of magnetic field
\begin{align}
    &  \int_0^{\tau} \int_{\Om} (b \p_t \phi+ b \vu \cdot \Grad \phi) \dxdt-\int_0^{\tau} \int_{ \Gamma^{ \rm{in} } }   b_{B}  \vu_{B}\cdot \mathbf{n} \phi \dS \dt      \nonumber     \\
    &   \quad  = \int_0^{\tau} \int_{ \Gamma^{ \rm{out} } }   b \vu_{B}\cdot \mathbf{n} \phi \dS \dt + \left[  \int_{\Om} b \phi \dx  \right]_{t=0}^{t=\tau}    \label{equ-magnetic}
\end{align}
      for any $\tau \in [0,T]$ and any $\phi\in C^{1}( [0,T] \times \overline{\Om})$;
}
\item{The energy inequality
\begin{align}
    &    \int_{\Om}    \left( \f{1}{2} \vr | \vu-\vu_{B}|^2 +\f{1}{\gamma-1} \vr^{\gamma} +\f{1}{2} b^2    \right) (\tau,\cdot) \dx   +  \int_0^{\tau} \int_{\Om}\mathbb{S}(\Grad \vu):\Grad \vu \dxdt    \nonumber     \\
    &   \quad   \quad        +  \int_0^{\tau} \int_{ \Gamma^{ \rm{out} } }  \left(  \f{1}{\gamma-1} \vr^{\gamma} +\f{1}{2} b^2  \right) \vu_{B} \cdot \mathbf{n}  \dS \dt    \nonumber     \\
    &   \quad     \leq      \int_{\Om}    \left( \f{1}{2} \vr_0 | \vu_0-\vu_{B}|^2 +\f{1}{\gamma-1} \vr_0^{\gamma} +\f{1}{2} b_0^2    \right)  \dx     \label{energy-inequa}           \\
    & \quad  \quad           -\int_0^{\tau} \int_{ \Gamma^{ \rm{in} } }  \left(  \f{1}{\gamma-1} \vr_{B}^{\gamma} +\f{1}{2} b_{B}^2  \right) \vu_{B} \cdot \mathbf{n}  \dS \dt     \nonumber    \\
    & \quad  \quad + \int_0^{\tau} \int_{\Om}  \left[   -\left( \vr^{\gamma}+\f{1}{2}b^2  \right) \Div \vu_{B} -\vr \vu \cdot \Grad \vu_{B} \cdot (\vu-\vu_{B}) +  \mathbb{S}(\Grad \vu):\Grad \vu_{B}                    \right] \dxdt   \nonumber
\end{align}
for a.e. $\tau\in (0,T)$.
}
\end{itemize}
\end{Definition}

Our first result concerning the existence of global-in-time weak solution to compressible viscous non-resistive MHD equations with general boundary data can be stated as
\begin{Theorem}\label{TH1}
Let $\Om\subset \R^2$ be a piecewise regular bounded Lipschitz domain and $\gamma>1$. Suppose that
\begin{align*}
    &   0< \vr_{B},b_{B} \in C_{c}(\R^2),\,\,  \vu_{B} \in C_{c}^1(\R^2;\R^2) ,       \\
    &    \int_{\Om}    \left( \f{1}{2} \vr_0 | \vu_0|^2 +\f{1}{\gamma-1} \vr_0^{\gamma} +\f{1}{2} b_0^2    \right)  \dx <\infty,    \\
    &   C_{\ast} \vr_0 \leq b_0 \leq C^{\ast} \vr_0,\, C_{\ast} \vr_{B} \leq b_{B} \leq C^{\ast} \vr_{B}   \,\,\text{    for some } 0<C_{\ast}<C^{\ast}<\infty .
\end{align*}

Then there exists a global weak solution to (\ref{in2})-(\ref{in4}) in the sense of Definition \ref{def:1}.
\end{Theorem}
\begin{Remark}
\begin{itemize}

\item{
In our previous work \cite{LS19}, existence of global weak solutions was proved for the same equations (\ref{in2}) with homogeneous Dirichlet boundary conditions. Therefore, Theorem \ref{TH1} extends
our previous result to the case of general inflow-outflow boundary conditions in a non-trivial manner. 
}

\item{
In \cite{LS19}, we are able to show the global well-posedness for the first level approximation equations, while this is not clear in the current context due to the general inflow-outflow boundary conditions, see Section \ref{con-rem} for more discussions. This is of independent interest.
}

\end{itemize}

\end{Remark}

Our second result concerns the stability of strong solutions in the class of weak solutions.
\begin{Theorem}\label{TH2}
Let $(\vr,\vu,b)$ be a weak solution to the problem (\ref{in2})-(\ref{in4}) ensured by Theorem \ref{TH1}. Assume that $(\tilde{\vr},\tilde{\vu},\tilde{b})$ is a strong solution to the same equations such that
\begin{align*}
    &   (\tilde{\vr},\tilde{\vu},\tilde{b})|_{t=0}=(\vr_0,\vu_0,b_0),         \\
    &   \tilde{\vu} |_{\p \Omega}= \vu_{B},\,\, \tilde{\vr} |_{ \Gamma^{ \text{in} } } =\vr_{B},\,\, \tilde{b} |_{ \Gamma^{ \text{in} } } =b_{B},    \\
    &   0<\underline{c}\leq \tilde{\vr},\tilde{b} \leq \overline{c}<\infty,   \\
    &   \tilde{\vr},\tilde{b} \in C^1 ([0,T]\times \overline{\Om}),\,\, \tilde{\vu} \in C^2 ([0,T]\times \overline{\Om};\R^2).
\end{align*}

Then
\[
\vr= \tilde{\vr},\,    \vu=\tilde{\vu}          , \,b=\tilde{b} \,\,\,\, \text{    a.e. in    } (0,T) \times \Om.
\]

\end{Theorem}

\begin{Remark}
\begin{itemize}
\item{
If the initial data and boundary data are smooth and obey the compatibility conditions, the problem (\ref{in2})-(\ref{in4}) admits local-in-time smooth solution. This may be proved by adapting the arguments from Valli and Zajaczkowski \cite{VaZa86}.
}
\item{
Following the same line of Theorem \ref{TH2} with slight modifications, one may establish the weak-strong uniqueness property for the compressible non-resistive MHD system (\ref{in2}) with homogeneous Dirichlet boundary conditions. This is not established in our previous work \cite{LS19}.
}
\end{itemize}
\end{Remark}

The rest of this paper is arranged as follows. In Section \ref{the-appro-eq}, we introduce the approximation scheme and discuss its global solvability. We then pass to the limits for approximating solutions in Section \ref{pass-limit}. This is the key to obtain the existence of global weak solutions. Remark that the specific mathematical structure of our model plays a crucial role in the compactness argument. In Section \ref{weak-strong}, we establish the relative energy inequality and prove the weak-strong uniqueness principle. Finally, we summarize the main results and give several comments or extensions, see Section \ref{con-rem}.

\section{Approximating equations}\label{the-appro-eq}

\subsection{The approximating equations}
Following Chang et al. \cite{CBN19}, we shall approximate the equation of continuity with artificial diffusion and Robin-type boundary conditions:
\begin{equation}\label{appro-conti}
\left\{\begin{aligned}
&  \p_t \vr +\Div (\vr \vu) =\ep \Delta \vr, \\
& \vr|_{t=0}=\vr_0,  \\
&  \ep \Grad \vr \cdot \mathbf{n} +(\vr_{B} -\vr ) [\vu_{B}\cdot \mathbf{n} ]^{-}|_{\p \Om}=0.    \\
\end{aligned}\right.
\end{equation}
Here and in what follows, we denote by $[f]^{-}:=\min\{ f,0 \}$. In a similar manner, we shall approximate the non-resistive magnetic equation through
\begin{equation}\label{appro-magn}
\left\{\begin{aligned}
&  \p_t b +\Div (b \vu) =\ep \Delta b, \\
& b|_{t=0}=b_0,  \\
&  \ep \Grad b \cdot \mathbf{n} +(b_{B} -b ) [\vu_{B}\cdot \mathbf{n} ]^{-}|_{\p \Om}=0.    \\
\end{aligned}\right.
\end{equation}
We shall use the Galerkin approximation for the balance of momentum. To this end, we first introduce a sequence of finite-dimensional spaces
\[
X_{n}=\text{span}\{  \mathbf{w}_i \, |\, \mathbf{w}_i \in C_c^{\infty}(\Om;\R^2) ,i=1,\ldots,n    \}
\]
and furthermore assume that $\{ \mathbf{w}_i \}_{i=1}^{\infty}$ are orthonormal with respect to the scalar product in $L^2(\Om;\R^2)$. The approximate velocity field takes the form
\[
\vu=\vv+\vu_{B},\,\, \vv \in C([0,T];X_n)
\]
and the approximate momentum equation reads
\begin{align}
    &   \int_0^{\tau}\int_{\Om} \left(  \vr\vu \cdot \p_t \vc{\phi} +\vr\vu \otimes \vu:\Grad \vc{\phi}+ \left( \vr^{\gamma}+\f{1}{2}b^2+       \delta (\vr+b)^{\beta} \right)   \Div \vc{\phi}                \right) \dxdt   \nonumber      \\
    &   \quad  -  \int_0^{\tau}\int_{\Om} \left(     \mathbb{S}(\Grad \vu) :\Grad \vc{\phi}+ \ep \Grad \vr \cdot \Grad \vu \cdot \vc{\phi}          \right) \dxdt
    =\left[  \int_{\Om}   \vr \vu \cdot  \vc{\phi}   \dx  \right]_{t=0}^{t=\tau}      \label{appro-mome}
\end{align}
for any $\vc{\phi}\in C^1([0,T];X_n)   $ with
\[
\vr \vu |_{t=0}=\vr_0 \vu_0,\,\, \vu_0=\vv_0+\vu_{B},\,\,\vv_0 \in X_n.
\]
Here, $\beta>1$ is a suitably large parameter. In particular, the exact value of $\beta$ has no influence in the limit processes and the final result.

However, due to the rough smoothness of the boundary, we have to use the approximate equation of continuity and magnetic field in the weak form:
\begin{align}
    &  \int_0^{\tau} \int_{\Om} (\vr \p_t \phi+ \vr \vu \cdot \Grad \phi  -\ep \Grad \vr \cdot \Grad \phi) \dxdt   - \int_0^{\tau} \int_{\p \Om} \phi \vr \vu_{B} \cdot \mathbf{n}   \dS \dt    \nonumber   \\
    &   \quad  +\int_0^{\tau} \int_{\p \Om} \phi   (\vr-\vr_{B})   [\vu_{B} \cdot \mathbf{n}]^{-}   \dS \dt = \left[  \int_{\Om} \vr \phi \dx  \right]_{t=0}^{t=\tau}  ,    \label{appro-conti-1}  \\
    & \vr|_{t=0}=\vr_0,      \nonumber
\end{align}
and
\begin{align}
    &  \int_0^{\tau} \int_{\Om} (b \p_t \phi+ b \vu \cdot \Grad \phi  -\ep \Grad b \cdot \Grad \phi) \dxdt   - \int_0^{\tau} \int_{\p \Om} \phi b \vu_{B} \cdot \mathbf{n}   \dS \dt    \nonumber   \\
    &   \quad  +\int_0^{\tau} \int_{\p \Om} \phi   (b-b_{B})   [\vu_{B} \cdot \mathbf{n}]^{-}   \dS \dt = \left[  \int_{\Om} b \phi \dx  \right]_{t=0}^{t=\tau}  ,    \label{appro-magn-1}  \\
    & b|_{t=0}=b_0,      \nonumber
\end{align}
for any $\tau \in [0,T]$ and any $\phi$ in the class
\[
\phi \in L^2(0,T;W^{1,2}(\Om)),\,\, \p_t \phi \in L^1(0,T;L^2(\Om)).
\]

\subsection{Solvability of the approximating equations}
Global solvability of the approximating problem (\ref{appro-mome})-(\ref{appro-magn-1}) is essentially based on the classical fixed point argument. To this end, we first fix the velocity field and
consider the solvability of the approximate equations of continuity and magnetic field. More precisely,
\begin{Lemma}\label{lem-3.1}
Let $\Om \subset \R^2$ be a bounded Lipschitz domain and $\vu=\vv+\vu_{B}$ with $ \vv \in C([0,T];X_n)$.
\begin{align*}
    &   0< \vr_{B},b_{B} \in C_{c}(\R^2),\,\,  \vu_{B} \in C_{c}^1(\R^2;\R^2) ,         \\
    &   0<\vr_0,b_0 \in C^1(\overline{\Om}) ,\,\,  \vu_0=\vv_0+\vu_{B},\,\,\vv_0\in C^1(\overline{\Om};\R^2).
\end{align*}
Then
\begin{itemize}
\item{
The approximate equation of continuity (\ref{appro-conti}) admits a weak solution $\vr$ in the sense of (\ref{appro-conti-1}), unique in the class
\[
\vr  \in C([0,T];L^2(\Om)) \cap    L^2(0,T;W^{1,2}(\Om))
\]
with the norm depending only on $\vr_0,\vr_{B},\vu_{B},\vv$. Besides, $\p_t \vr\in L^2((0,T)\times \Om),\sqrt{\ep} \Grad \vr \in L^{\infty}(0,T;L^2(\Om))$ with the norm depending only on $\vr_0,\vr_{B},\vu_{B},\vv$.         \\
The approximate equation of magnetic field (\ref{appro-magn}) admits a weak solution $b$ in the sense of (\ref{appro-magn-1}), unique in the class
\[
b  \in C([0,T];L^2(\Om)) \cap    L^2(0,T;W^{1,2}(\Om))
\]
with the norm depending only on $b_0,b_{B},\vu_{B},\vv$. Besides, $\p_t b \in L^2((0,T)\times \Om),\sqrt{\ep} \Grad b \in L^{\infty}(0,T;L^2(\Om))$ with the norm depending only on $b_0,b_{B},\vu_{B},\vv$.
}
\item{
For any $\tau \in [0,T]$,
\begin{align*}
    & \| \vr(\tau) \|_{L^{\infty}(\Om)} \leq M_{\vr}\exp (T \| \Div \vu \|_{ L^{\infty}( (0,\tau)\times \Om )    }  ),           \\
     & \| \vr^{-1}(\tau) \|_{L^{\infty}(\Om)} \leq m_{\vr}^{-1}   \exp (T \| \Div \vu \|_{ L^{\infty}( (0,\tau)\times \Om )    }  ),           \\
    &    \| b(\tau) \|_{L^{\infty}(\Om)} \leq M_{b}\exp (T \| \Div \vu \|_{ L^{\infty}( (0,\tau)\times \Om )    }  ), \\
     & \| b^{-1}(\tau) \|_{L^{\infty}(\Om)} \leq m_{b}^{-1}   \exp (T \| \Div \vu \|_{ L^{\infty}( (0,\tau)\times \Om )    }  ).
\end{align*}
For a.e. $\tau \in (0,T)$ and a.e. $x\in \p \Om$,
\begin{align*}
    & \vr(\tau,x) \leq M_{\vr}\exp (T \| \Div \vu \|_{ L^{\infty}( (0,\tau)\times \Om )    }  ),           \\
    &  \vr^{-1}(\tau,x) \leq         m_{\vr}^{-1}   \exp (T \| \Div \vu \|_{ L^{\infty}( (0,\tau)\times \Om )    }  ),                         \\
    &    b(\tau,x) \leq M_{b}\exp (T \| \Div \vu \|_{ L^{\infty}( (0,\tau)\times \Om )    }  ),\\
    &  b^{-1}(\tau,x) \leq         m_{b}^{-1}   \exp (T \| \Div \vu \|_{ L^{\infty}( (0,\tau)\times \Om )    }  ).
\end{align*}
Here
\begin{align*}
    &    M_{\vr}:= \max\{  \max_{\Om}\vr_0, \max_{\Gamma^{\rm in}} \vr_{B},  \| \vu_{B} \|_{ L^{\infty}( (0,T) \times \Om )    }             \} ,       \\
    &   M_{b}:= \max\{  \max_{\Om} b_0, \max_{\Gamma^{\rm in}} b_{B},  \| \vu_{B} \|_{ L^{\infty}( (0,T) \times \Om )    }             \} , \\
    &     m_{\vr}:= \min\{  \min_{\Om}\vr_0, \min_{\Gamma^{\rm in}} \vr_{B}  \}  ,          \\
       &     m_{b}:= \min\{  \min_{\Om}b_0, \min_{\Gamma^{\rm in}} b_{B}  \}  .
\end{align*}
}
\item{
The weak formulations (\ref{appro-conti-1})-(\ref{appro-magn-1}) are satisfied in the renormalized sense, i.e.,
\begin{align*}
    &   \ep \int_0^{\tau} \int_{\Om}   \Phi''(\vr) |\Grad \vr|^2  \dxdt + \int_0^{\tau} \int_{\Om} \Div(\vr \vu) \Phi'(\vr)    \dxdt      \\
    &   \quad  + \left[  \int_{\Om} \Phi (\vr)\dx \right]_{t=0}^{t=\tau} =   \int_0^{\tau} \int_{\Om}          \Phi'(\vr)(\vr-\vr_{B}) [\vu_{B} \cdot \mathbf{n}]^{-}   \dS \dt,
\end{align*}
\begin{align*}
    &   \ep \int_0^{\tau} \int_{\Om}   \Phi''(b) |\Grad b|^2  \dxdt + \int_0^{\tau} \int_{\Om} \Div(b \vu) \Phi'(b)    \dxdt      \\
    &   \quad  + \left[  \int_{\Om} \Phi (b)\dx \right]_{t=0}^{t=\tau} =   \int_0^{\tau} \int_{\Om}          \Phi'(b)(b-b_{B}) [\vu_{B} \cdot \mathbf{n}]^{-}   \dS \dt
\end{align*}
for any $\tau \in [0,T]$ and any $\Phi \in C^2(\R)$.
}
\end{itemize}
\end{Lemma}

The proof of Lemma \ref{lem-3.1} follows the same line as section 3.1 of Abbatiello et al. \cite{AFN21} or section 4.1 of Kra\v{c}mar et al. \cite{KKNN22}. For simplicity, we omit the details here.

Inspired by the compressible Navier-Stokes equations with general inflow-outflow boundary conditions \cite{CBN19} and the compressible multifluid model with general inflow-outflow boundary conditions \cite{KKNN22},
we first solve the approximate equation of continuity (\ref{appro-conti-1}) and magnetic field (\ref{appro-magn-1}), for given $\vu=\vv+\vu_{B},\,\vv\in C([0,T];X_n)$, and obtain the unique solution respectively as $\vr=\vr[\vu],b=b[\vu]$; the later is then inserted into the approximate equation of momentum (\ref{appro-mome}) in order to obtain a unique solution $\vu=\vu[ \vr[\vu],b[\vu]]$. This indeed defines a mapping from
$C([0,T];X_n)$ to itself via
\[
\vv \mapsto \vu[ \vr[\vu],b[\vu]]- \vu_{B}.
\]
The classical fixed point argument is then applied to the mapping to conclude that (\ref{appro-mome})-(\ref{appro-magn-1}) are globally solvable. To summarize, we have
\begin{Proposition}\label{prop-1}
Let $\Om \subset \R^2$ be a bounded Lipschitz domain and
\begin{align*}
    &   0< \vr_{B},b_{B} \in C_{c}(\R^2),\,\,  \vu_{B} \in C_{c}^1(\R^2;\R^2) ,         \\
    &   0<\vr_0,b_0 \in C^1(\overline{\Om}) .
\end{align*}
Then for any fixed $\ep,\delta>0,n\geq 1$, there exists a solution $(\vr_n,\vu_n,b_n)$, with $\vu_n=\vv_n+\vu_{B},\vv_n \in C([0,T];X_n)$, to the approximate problem (\ref{appro-conti})-(\ref{appro-magn-1}) in
the time-space cylinder $(0,T)\times \Om$. Moreover,
\begin{itemize}
\item{
Lower and upper bounds for the density and magnetic field:
\begin{align}
    &   \vr_n(t,x)\geq c(n)>0,\,\,  b_n(t,x)\geq c(n)>0,   \nonumber     \\
    &   C_{\ast} \vr_n(t,x) \leq b_n(t,x)\leq C^{\ast}\vr_n(t,x)     \label{sov-1}
\end{align}
for any $t \in [0,T]$ and a.e. $x\in \Om$;
\begin{align}
    &   \vr_n(t,x)\geq c(n)>0,\,\,  b_n(t,x)\geq c(n)>0,    \nonumber      \\
    &   C_{\ast} \vr_n(t,x) \leq b_n(t,x)\leq C^{\ast}\vr_n(t,x)    \label{sov-2}
\end{align}
for a.e. $t \in (0,T)$ and a.e. $x\in \p \Om$.
}
\item{
The approximate energy inequality
\begin{align}
    &  \left[   \int_{\Om} \left(   \f{1}{2}\vr_n |\vv_n|^2 +\f{ \vr_n^{\gamma} }{\gamma-1} +\f{1}{2}b_n^2+\f{\delta}{\beta-1}(\vr_n+b_n)^{\beta}       \right) \dx      \right]_{t=0}^{t=\tau}
    + \int_0^{\tau} \int_{\Om} \mathbb{S}(\Grad \vu_n):\Grad \vu_n \dxdt    \nonumber       \\
    &   \quad  \quad   +  \int_0^{\tau} \int_{\Om}  \Big(   \ep \gamma \vr_n ^{\gamma-2}|\Grad \vr_n|^2 +\ep |\Grad b_n|^2 +\ep \delta \beta (\vr_n+b_n)^{\beta-2} |\Grad(\vr_n+b_n)|^2 \Big) \dxdt   \nonumber    \\
    &   \quad \quad  +    \int_0^{\tau} \int_{ \Gamma^{ \rm{in} } }  \f{1}{\gamma-1}\Big( \vr_{B}^{\gamma}-\gamma\vr_{n}^{\gamma-1}(\vr_{B}-\vr_n ) -\vr_n^{\gamma}       \Big)  |\vu_{B} \cdot \mathbf{n}|  \dS\dt   \nonumber                    \\
    &   \quad \quad  +    \int_0^{\tau} \int_{ \Gamma^{ \rm{in} } }       \Big(   \f{1}{2}b_{B}^2 -b_n(b_{B}-b_n)-\f{1}{2}b_n^2               \Big)  |\vu_{B} \cdot \mathbf{n}|  \dS\dt  \nonumber      \\
    &  \quad \quad  +    \int_0^{\tau} \int_{ \Gamma^{ \rm{in} } }  \f{\delta}{\beta-1}     \Big( (\vr_B +b_B)^{\beta} -\beta (\vr_n+b_n)^{\beta-1} (\vr_B+b_B-\vr_n-b_n)   -  (\vr_n+b_n)^{\beta}   \Big)  |\vu_{B} \cdot \mathbf{n}|  \dS\dt             \nonumber            \\
    & \quad   \quad  +  \int_0^{\tau} \int_{ \Gamma^{ \rm{out} } }      \left( \f{ \vr_n^{\gamma} }{\gamma-1} +\f{1}{2}b_n^2+\f{\delta}{\beta-1}(\vr_n+b_n)^{\beta}       \right) \vu_{B} \cdot \mathbf{n}   \dS\dt \nonumber     \\
    & \quad   \leq    \int_0^{\tau} \int_{\Om} \mathbb{S}(\Grad \vu_n):\Grad \vu_B \dxdt   -   \int_0^{\tau} \int_{ \Gamma^{ \rm{in} } }     \left( \f{ \vr_B^{\gamma} }{\gamma-1} +\f{1}{2}b_B^2+\f{\delta}{\beta-1}(\vr_B+b_B)^{\beta}       \right) \vu_{B} \cdot \mathbf{n}   \dS\dt             \nonumber            \\
    &      \quad   \quad \quad -    \int_0^{\tau} \int_{\Om}     \Big[ \vr_n \vu_n \otimes \vu_n +\Big(  \vr_n^{\gamma}+\f{1}{2}b_n^2+\delta(\vr_n+b_n)^{\beta}   \Big)\mathbb{I}             \Big]:\Grad \vu_B \dxdt  \nonumber     \\
    &    \quad   \quad \quad   +\f{1}{2}\int_0^{\tau} \int_{\Om}  \vr_n \vu_n \cdot \Grad |\vu_B|^2 \dxdt +\int_0^{\tau} \int_{\Om} \ep \Grad \vr_n \cdot \Grad \vv_n \cdot \vu_B \dxdt           \label{sov-3}
\end{align}
holds for any $\tau\in [0,T]$.
}
\end{itemize}

\end{Proposition}
{\emph {Proof of Proposition \ref{prop-1}.}}  Observe that the lower bounds for the density and magnetic field directly follow from the second item of Lemma \ref{lem-3.1}. Furthermore, the domination relations between the density and magnetic field hold in the same spirit upon seeing that $b_n-C_{\ast}\vr_n$ and $C^{\ast}\vr_n -b_n$ satisfy the same type of equations.
We are thus concentrate on the derivation of approximate energy inequality. It follows from Lemma \ref{lem-3.1} that
\begin{align}
    &   \left[ \int_{\Om} \f{\vr_n^{\gamma}}{\gamma-1} \dx\right]_{t=0}^{t=\tau} +\int_0^{\tau} \int_{\Om}    \ep \gamma \vr_n ^{\gamma-2}|\Grad \vr_n|^2 \dxdt
    +  \int_0^{\tau} \int_{ \Gamma^{ \rm{out} } }      \f{\vr_n^{\gamma}}{\gamma-1}           \vu_{B} \cdot \mathbf{n}   \dS\dt     \nonumber    \\
    &   \quad   \quad +  \int_0^{\tau} \int_{ \Gamma^{ \rm{in} } }    \f{1}{\gamma-1}\Big( \vr_{B}^{\gamma}-\gamma\vr_{n}^{\gamma-1}(\vr_{B}-\vr_n ) -\vr_n^{\gamma}       \Big)  |\vu_{B} \cdot \mathbf{n}|  \dS\dt    \nonumber    \\
    &   \quad    =  - \int_0^{\tau} \int_{\Om} \vr_n^{\gamma} \Div \vu_n \dxdt -   \int_0^{\tau} \int_{ \Gamma^{ \rm{in} } }     \f{ \vr_B^{\gamma} }{\gamma-1}    \vu_{B} \cdot \mathbf{n}      \dS\dt ,  \label{sov-4}
\end{align}
\begin{align}
    &   \left[ \int_{\Om} \f{1}{2}   b_n^{2} \dx\right]_{t=0}^{t=\tau} +\int_0^{\tau} \int_{\Om}    \ep |\Grad b_n|^2 \dxdt
    +  \int_0^{\tau} \int_{ \Gamma^{ \rm{out} } }     \f{1}{2}   b_n^{2}         \vu_{B} \cdot \mathbf{n}   \dS\dt     \nonumber    \\
    &   \quad   \quad +  \int_0^{\tau} \int_{ \Gamma^{ \rm{in} } }       \Big(   \f{1}{2}b_{B}^2 -b_n(b_{B}-b_n)-\f{1}{2}b_n^2               \Big)  |\vu_{B} \cdot \mathbf{n}|  \dS\dt    \nonumber    \\
    &   \quad    =  - \int_0^{\tau} \int_{\Om} \f{1}{2}   b_n^{2} \Div \vu_n \dxdt -   \int_0^{\tau} \int_{ \Gamma^{ \rm{in} } }     \f{1}{2} b_B^{2}     \vu_{B} \cdot \mathbf{n}      \dS\dt ,  \label{sov-5}
\end{align}
\begin{align}
    &   \left[ \int_{\Om} \f{(\vr_n+b_n)^{\beta}}{\beta-1} \dx\right]_{t=0}^{t=\tau} +\int_0^{\tau} \int_{\Om}   \ep \delta \beta (\vr_n+b_n)^{\beta-2} |\Grad(\vr_n+b_n)|^2  \dxdt
       \nonumber    \\
    &   \quad   \quad +  \int_0^{\tau} \int_{ \Gamma^{ \rm{in} } }  \f{\delta}{\beta-1}     \Big( (\vr_B +b_B)^{\beta} -\beta (\vr_n+b_n)^{\beta-1} (\vr_B+b_B-\vr_n-b_n)   -  (\vr_n+b_n)^{\beta}   \Big)  |\vu_{B} \cdot \mathbf{n}|  \dS\dt       \nonumber    \\
    &     \quad   \quad    +  \int_0^{\tau} \int_{ \Gamma^{ \rm{out} } }     \f{\delta}{\beta-1}(\vr_n+b_n)^{\beta}          \vu_{B} \cdot \mathbf{n}   \dS\dt          \nonumber    \\
    &   \quad    =  - \int_0^{\tau} \int_{\Om}      \delta (\vr_n +b_n)^{\beta}   \Div \vu_n \dxdt -   \int_0^{\tau} \int_{ \Gamma^{ \rm{in} } }    \f{\delta}{\beta-1}(\vr_B+b_B)^{\beta}      \vu_{B} \cdot \mathbf{n}      \dS\dt .\label{sov-6}
\end{align}
Next, we test the approximate momentum equation (\ref{appro-mome}) by $\vv_n$ to obtain
\begin{align}
    &   \int_0^{\tau} \int_{\Om}  \Big(  \p_t (\vr_n \vv_n) \cdot \vv_n-\vr_n \vu_n \otimes \vu_n :\Grad\vv_n         \Big) \dxdt   +    \int_0^{\tau} \int_{\Om}  \p_t \vr_n \vu_B \cdot \vv_n      \dxdt    \nonumber  \\
    &   \quad  \quad  + \int_0^{\tau} \int_{\Om}    \Big(   \mathbb{S}(\Grad \vu_n):\Grad \vv_n +\ep \Grad\vr_n \cdot \Grad \vu_n \cdot \vv_n    \Big) \dxdt        \nonumber    \\
    &   \quad    =  \int_0^{\tau} \int_{\Om}    \left( \vr_n^{\gamma}+\f{1}{2}b_n^2+       \delta (\vr_n+b_n)^{\beta} \right)    \Div \vv_n    \dxdt .        \label{sov-7}
\end{align}

The approximate energy inequality then follows from (\ref{sov-4})-(\ref{sov-7}) as long as we invoke the basic calculations below:
\begin{align*}
    &        \int_0^{\tau} \int_{\Om}    \Big(  \p_t (\vr_n \vv_n) \cdot \vv_n-\vr_n \vu_n \otimes \vu_n :\Grad\vv_n         \Big) \dxdt          \\
    &   \quad  =    \int_0^{\tau} \int_{\Om}    \Big(   \f{1}{2} \p_t(\vr_n |\vv_n|^2) -\ep \Grad \vr_n \cdot \Grad \vv_n \cdot \vv_n - \vr_n \vu_n \cdot  \Grad \vv_n \cdot \vu_B  \Big) \dxdt,
\end{align*}
\begin{align*}
    \int_0^{\tau} \int_{\Om}  \p_t \vr_n \vu_B \cdot \vv_n      \dxdt &    = \int_0^{\tau} \int_{\Om}   \Big(  \vr_n \vu_n \cdot \Grad \vv_n \cdot \vu_B+  \vr_n \vu_n \cdot \Grad \vu_B \cdot \vv_n                   \\
    &   \quad  \quad \quad     - \ep \Grad \vr_n \cdot \Grad \vu_B \cdot \vv_n-\ep \Grad \vr_n \cdot \Grad \vv_n \cdot \vu_B \Big)   \dxdt  ,
\end{align*}
where we have conveniently used the fact that $\vu_n=\vv_n+\vu_{B}$ and integration by parts. This completes the proof of Proposition \ref{prop-1}.   $\Box$

\section{Passing to the limits}\label{pass-limit}
To obtain the existence of global weak solutions, we shall successively pass to the limits for $n\rightarrow \infty,\ep\rightarrow 0$ and $\delta \rightarrow 0$. Indeed, the limit process
$n\rightarrow \infty$ is more or less simpler than the latter ones thanks to the nice regularity at this stage, which could be carried out as \cite{FNP,KKNN22}. The limit process $\ep\rightarrow 0$ and that of $\delta \rightarrow 0$ share certain similar features.
It is also known that the last step is more delicate due to the loss of higher integrability of density and magnetic field. Therefore, we shall concentrate on the limit process $\delta \rightarrow 0$ to simplify the
presentation.

In accordance with the previous discussions, after passing to the limits for $n\rightarrow \infty$ and $\ep\rightarrow 0$ we may assume that there exists a sequence of approximate solutions $\{\vr_{\delta},\vu_{\delta},b_{\delta} \}_{\delta>0}$ obeying
\begin{itemize}
\item{
\begin{align}
    &    \vr_{\delta},b_{\delta} \in C_{\rm weak}([0,T];L^{\beta}(\Om)  ) \cap L^{\beta}( 0,T;L^{\beta} (\Gamma^{\rm out};  \vu_{B} \cdot \mathbf{n} \dS )    ),     \nonumber       \\
    &  \vv_{\delta}:= \vu_{\delta}-\vu_{B} \in L^2(0,T; W_0^{1,2} (\Om;\R^2)),\,\, \vr_{\delta} \vu_{\delta} \in C_{\rm weak}([0,T];L^{  \f{2\beta}{\beta+1} }(\Om)  ),  \label{pass-1}   \\
    &   \vr_{\delta} |\vv_{\delta}|^2 \in L^{\infty}(0,T; L^1(\Om))  ;  \nonumber
\end{align}
}
\item{
\begin{align}
    &   \vr_{\delta}(t,x),b_{\delta}(t,x) \geq 0,   \nonumber     \\
    &   C_{\ast} \vr_{\delta}(t,x) \leq b_{\delta}(t,x)\leq C^{\ast}\vr_{\delta}(t,x)  \,\, \text{  for any }t\in [0,T], \text{ a.e. }x\in \Om,  \label{pass-2}  \\
     &   \vr_{\delta}(t,x),b_{\delta}(t,x) \geq 0,   \nonumber     \\
    &   C_{\ast} \vr_{\delta}(t,x) \leq b_{\delta}(t,x)\leq C^{\ast}\vr_{\delta}(t,x)  \,\, \text{  for a.e. }t\in (0,T), \text{ a.e. }x\in \p \Om;  \nonumber
\end{align}
}
\item{
\begin{align}
    &  \int_0^{\tau} \int_{\Om} \Big(\vr_{\delta} \p_t \phi+ \vr_{\delta} \vu_{\delta} \cdot \Grad \phi  \Big) \dxdt-\int_0^{\tau} \int_{ \Gamma^{ \rm{in} } }   \vr_{B}  \vu_{B}\cdot \mathbf{n} \phi \dS \dt    \nonumber      \\
    &   \quad  = \int_0^{\tau} \int_{ \Gamma^{ \rm{out} } }   \vr_{\delta}  \vu_{B}\cdot \mathbf{n} \phi \dS \dt + \left[  \int_{\Om} \vr_{\delta} \phi \dx  \right]_{t=0}^{t=\tau}  \label{pass-3}
\end{align}
      for any $\tau \in [0,T]$ and any $\phi\in C^{1}( [0,T] \times \overline{\Om})$;
}
\item{
\begin{align}
    &  \int_0^{\tau} \int_{\Om} \Big(   \vr_{\delta} \vu_{\delta} \cdot \p_t \vc{\phi} +\vr_{\delta} \vu_{\delta} \otimes \vu_{\delta} :\Grad \vc{\phi}+\left( \vr_{\delta}^{\gamma}+\f{1}{2}b_{\delta}^2
     +\delta (\vr_{\delta}+b_{\delta})^{\beta}    \right)\Div \vc{\phi}     \nonumber    \\
    &   \quad \quad     -\mathbb{S}(\Grad \vu_{\delta} ):\Grad \vc{\phi}                     \Big) \dxdt  =  \left[  \int_{\Om} \vr_{\delta} \vu_{\delta} \cdot \vc{\phi} \dx  \right]_{t=0}^{t=\tau}    \label{pass-4}
\end{align}
      for any $\tau \in [0,T]$ and any $\vc{\phi} \in C_c^{1}( [0,T] \times \Om;\R^2)$;
}
\item{
\begin{align}
    &  \int_0^{\tau} \int_{\Om} \Big(b_{\delta} \p_t \phi+ b_{\delta} \vu_{\delta} \cdot \Grad \phi \Big) \dxdt-\int_0^{\tau} \int_{ \Gamma^{ \rm{in} } }   b_{B}  \vu_{B}\cdot \mathbf{n} \phi \dS \dt   \nonumber       \\
    &   \quad  = \int_0^{\tau} \int_{ \Gamma^{ \rm{out} } }   b_{\delta} \vu_{B}\cdot \mathbf{n} \phi \dS \dt + \left[  \int_{\Om} b_{\delta} \phi \dx  \right]_{t=0}^{t=\tau}  \label{pass-5}
\end{align}
      for any $\tau \in [0,T]$ and any $\phi\in C^{1}( [0,T] \times \overline{\Om})$;
}
\item{
\begin{align}
    &  \left[   \int_{\Om} \left(   \f{1}{2}\vr_{\delta} |\vv_{\delta}|^2 +\f{ \vr_{\delta}^{\gamma} }{\gamma-1} +\f{1}{2}b_{\delta}^2+\f{\delta}{\beta-1}(\vr_{\delta}+b_{\delta})^{\beta}       \right) \dx      \right]_{t=0}^{t=\tau}
    + \int_0^{\tau} \int_{\Om} \mathbb{S}(\Grad \vu_{\delta}):\Grad \vu_{\delta} \dxdt    \nonumber       \\
    & \quad   \quad  +  \int_0^{\tau} \int_{ \Gamma^{ \rm{out} } }      \left( \f{ \vr_{\delta}^{\gamma} }{\gamma-1} +\f{1}{2}b_{\delta}^2+\f{\delta}{\beta-1}(\vr_{\delta}+b_{\delta})^{\beta}       \right) \vu_{B} \cdot \mathbf{n}   \dS\dt \nonumber     \\
    & \quad   \leq    \int_0^{\tau} \int_{\Om} \mathbb{S}(\Grad \vu_{\delta}):\Grad \vu_B \dxdt   -   \int_0^{\tau} \int_{ \Gamma^{ \rm{in} } }     \left( \f{ \vr_B^{\gamma} }{\gamma-1} +\f{1}{2}b_B^2+\f{\delta}{\beta-1}(\vr_B+b_B)^{\beta}       \right) \vu_{B} \cdot \mathbf{n}   \dS\dt             \nonumber            \\
    &      \quad   \quad \quad -    \int_0^{\tau} \int_{\Om}     \Big[ \vr_{\delta} \vu_{\delta} \otimes \vu_{\delta} +\Big(  \vr_{\delta}^{\gamma}+\f{1}{2}b_{\delta}^2+\delta(\vr_{\delta}+b_{\delta})^{\beta}   \Big)\mathbb{I}             \Big]:\Grad \vu_B \dxdt  \nonumber     \\
    &    \quad   \quad \quad   +\f{1}{2}\int_0^{\tau} \int_{\Om}  \vr_{\delta} \vu_{\delta} \cdot \Grad |\vu_B|^2 \dxdt            \label{pass-6}
\end{align}
for a.e. $\tau\in (0,T)$.
}

\end{itemize}

\subsection{Uniform-in-$\delta$ estimates}
It follows readily from the domination condition (\ref{pass-2}) and the energy inequality (\ref{pass-6}) that
\begin{align}
    &    \| \vu_{\delta} \|_{ L^2(0,T;W^{1,2}(\Om))   } \leq C,       \nonumber      \\
    &   \|   \vr_{\delta} |\vv_{\delta}|^2   \|_{  L^{\infty}(0,T;L^1(\Om))   } \leq C,   \nonumber    \\
    &   \|      (\vr_{\delta},b_{\delta})      \|_{  L^{\infty}(0,T;L^{  \max\{ 2,\gamma \}   }(\Om))      }  \leq C,      \nonumber         \\
     &  \delta \|      (\vr_{\delta},b_{\delta})      \|_{  L^{\infty}(0,T;L^{  \beta   }(\Om))      }^{\beta}  \leq C,         \label{pass-7}       \\
     &   \|          \vr_{\delta}|\vu_{B} \cdot \mathbf{n}|^{\f{1}{\gamma}}       \|_{ L^{\gamma}(0,T;L^{\gamma}(\p \Omega))          }  \leq C,       \nonumber          \\
     &     \|          b_{\delta}|\vu_{B} \cdot \mathbf{n}|^{\f{1}{2}}       \|_{ L^{2}(0,T;L^{2}(\p \Omega))          }  \leq C.         \nonumber
\end{align}
Here and in the sequel, we denote by $C$ a generic positive constant independent of $\delta$. Based on the uniform estimates (\ref{pass-7}), we conclude that, up to a suitable subsequence, there exists a weak limit $(\vr,\vu,b)$ such that as $\delta \rightarrow 0$
\begin{align}
    &   (\vr_{\delta},b_{\delta}) \rightarrow (\vr,b) \text{ in } C_{\rm weak}([0,T];L^{ \max\{ 2,\gamma \}  } (\Om)  ) ,    \nonumber       \\
    &  \vr_{\delta}  \rightarrow \vr     \text{ weakly in }  L^{\max\{ 2,\gamma \} }(0,T; L^{\max\{ 2,\gamma \}}(\p \Omega;|\vu_{B}\cdot \mathbf{n}|\dS)            ) ,      \nonumber               \\
  &  b_{\delta}  \rightarrow    b     \text{ weakly in }  L^{\max\{ 2,\gamma \} }(0,T; L^{\max\{ 2,\gamma \}}(\p \Omega;|\vu_{B}\cdot \mathbf{n}|\dS)            ) ,         \label{pass-8}            \\
    &  \vu_{\delta}       \rightarrow \vu   \text{ weakly in }  L^2(0,T;W^{1,2}(\Om;\R^2)) ,\,\, \vu-\vu_{B}=:\vv \in   L^2(0,T;W_{0}^{1,2}(\Om;\R^2)),     \nonumber         \\
    &     \vr_{\delta}   \vu_{\delta}       \rightarrow          \vr \vu \text{ in } C_{\rm weak}([0,T];L^{   \f{2\max\{ 2,\gamma \}}{\max\{ 2,\gamma \}+1}      } (\Om)  ) .   \nonumber  \\
    &        \vr_{\delta}   \vu_{\delta} \otimes \vu_{\delta}      \rightarrow     \vr \vu\otimes \vu    \text{ in }  \mathcal{D}'((0,T)\times \Om),   \nonumber        \\
    &        b_{\delta}   \vu_{\delta}       \rightarrow   b\vu             \text{ in }  \mathcal{D}'((0,T)\times \Om).   \nonumber
\end{align}

Unfortunately, the convergence results in (\ref{pass-8}) are not sufficient to pass to the limit in the weak formulations due to in particular the pressure term. To this end, we first need to improve the integrability of
density and magnetic field. This is the aim of the next subsection.

\subsection{Pressure estimates}
In contrast to our previous paper \cite{LS19}, where the homogeneous Dirichlet boundary condition was treated, the present setting is more delicate. This turns out to be possible only in compact subsets of $\Om$
due to the general inflow-outflow boundary conditions. More precisely, for any $0\leq \eta \in C_c^{\infty}((0,T))$ and $0\leq \psi \in C_c^{\infty}(\Om)$, we use\footnote{Here, we denote by $\mathscr{B}$ the well-known Bogovskii operator, which in particular solves
\[
\Div \mathscr{B}[f]=f  \text{  in  }\Om,\,\,\,\, \mathscr{B}[f]|_{\p \Om}=0,
\]
and could be extended as a bounded linear operator
\[
\mathscr{B}: \left\{  f\in L^p(\Om)\,:\, \int_{\Om}f \dx=0         \right\}\rightarrow W_0^{1,p}(\Om;\R^2),\,\,\,\,1<p<\infty.
\]
Moreover, $\|\mathscr{B}[f] \|_{W_0^{1,p}(\Om)   } \leq C \| f \|_{L^p(\Om)}  $ for any $1<p<\infty$. For more discussions and properties of Bogovskii operator, we refer to \cite{FNP}.
}
\[
\vc{\phi}(t,x)=\eta(t) \mathscr{B}\left[  \psi  \vr_{\delta}^{\vartheta}- \f{1}{|\Om|} \int_{\Om}   \psi \vr_{\delta}^{\vartheta}  \dx              \right],\,\,\,\, \vartheta>0 \text{  is sufficiently small }
\]
as a test function in the momentum equation (\ref{pass-4}) to arrive at
\begin{align}
    &  \int_0^T \eta \int_{\Om}      \left( \vr_{\delta}^{\gamma}+\f{1}{2}b_{\delta}^2
     +\delta (\vr_{\delta}+b_{\delta})^{\beta}    \right)        \vr_{\delta}^{\vt} \psi \dxdt    \nonumber        \\
    &   \quad  =  \int_0^T  \left[  \eta     \f{1}{|\Om|} \int_{\Om}   \psi \vr_{\delta}^{\vartheta}  \dx   \int_{\Om}      \left( \vr_{\delta}^{\gamma}+\f{1}{2}b_{\delta}^2
     +\delta (\vr_{\delta}+b_{\delta})^{\beta}    \right)        \dx                \right] \dt             \nonumber             \\
    &   \quad  \quad   -      \int_0^T \eta'(t)    \int_{\Om}  \vr_{\delta}\vu_{\delta} \cdot \mathscr{B}\left[  \psi  \vr_{\delta}^{\vartheta}- \f{1}{|\Om|} \int_{\Om}   \psi \vr_{\delta}^{\vartheta}  \dx              \right]  \dxdt                     \nonumber         \\
   &   \quad  \quad   -      \int_0^T \eta(t)    \int_{\Om}  \vr_{\delta}\vu_{\delta} \cdot \mathscr{B}\left[  \p_t \left( \psi  \vr_{\delta}^{\vartheta}- \f{1}{|\Om|} \int_{\Om}   \psi \vr_{\delta}^{\vartheta}  \dx  \right)            \right]  \dxdt          \label{pass-9}                    \\
    &   \quad  \quad   -      \int_0^T \eta    \int_{\Om}  \vr_{\delta}\vu_{\delta} \otimes \vu_{\delta}:\Grad \mathscr{B}\left[  \psi  \vr_{\delta}^{\vartheta}- \f{1}{|\Om|} \int_{\Om}   \psi \vr_{\delta}^{\vartheta}  \dx              \right]  \dxdt                     \nonumber         \\
    &   \quad  \quad   +      \int_0^T \eta    \int_{\Om}     \mathbb{S} (\Grad \vu_{\delta}) :\Grad \mathscr{B}\left[  \psi  \vr_{\delta}^{\vartheta}- \f{1}{|\Om|} \int_{\Om}   \psi \vr_{\delta}^{\vartheta}  \dx              \right]  \dxdt    .    \nonumber
\end{align}
Essentially based on the uniform estimates (\ref{pass-7}), we conclude, via Sobolev embeddings, H\"{o}lder inequality, integration by parts, the renormalized equations and the property of Bogovskii operator that the right-hand side of (\ref{pass-9}) is bounded
on compact subsets of $\Om$. This process is nowadays standard and we refer to \cite{FNP,KKNN22} for the details. Therefore, after a straightforward manipulation we obtain that there exists $\vt>0$ depending only on $\gamma$ such that
\beq\label{pass-10}
\int_0^T \int_{K} \left( \vr_{\delta}^{\gamma+\vt}+ b_{\delta}^{2+\vt} +\delta \vr_{\delta}^{\beta+\vt}   +\delta b_{\delta}^{\beta+\vt}          \right) \dxdt
\leq C(K)
\eeq
for any compact $K\subset \Om$. Here, we use $C(K)$ to emphasize the dependence on $K$ for generic positive constant. Thanks to the domination condition between the density and magnetic field (\ref{pass-2}), it holds further that
\beq\label{pass-11}
\int_0^T \int_{K}
\left( \vr_{\delta}^{   \max\{ \gamma+\vt,2+\vt \}    }  +b_{\delta}^{   \max\{ \gamma+\vt,2+\vt \}    }      \right) \dxdt \leq C(K)
\eeq
for any compact $K\subset \Om$. It follows from (\ref{pass-11}) that
\begin{align}
    &  \delta (\vr_{\delta}+b_{\delta})^{\beta}  \rightarrow 0    \text{ strongly in }  L^{1 }((0,T)\times K  ) ,     \nonumber        \\
    &  \vr_{\delta}^{\gamma}+\f{1}{2}b_{\delta}^2  \rightarrow  \overline{  \vr^{\gamma}+\f{1}{2} b^2    }   \text{ weakly in }  L^{      \min\left\{ \f{ \max\{ \gamma+\vt,2+\vt \} }{\gamma} ,\f{ \max\{ \gamma+\vt,2+\vt \} }{2}   \right\}            }((0,T)\times K  )          \label{pass-12}
\end{align}
for any compact $K\subset \Om$. Therefore, in combination with (\ref{pass-12}) we may pass to the limit $\delta \rightarrow 0$ in the approximate momentum equation (\ref{pass-4}) to obtain
\begin{align}
    &  \int_0^{\tau} \int_{\Om} \Big(   \vr \vu \cdot \p_t \vc{\phi} +\vr \vu \otimes \vu :\Grad \vc{\phi}+ \overline{  \vr^{\gamma}+\f{1}{2} b^2    } \,  \Div \vc{\phi}     \nonumber    \\
    &   \quad \quad     -\mathbb{S}(\Grad \vu ):\Grad \vc{\phi}                     \Big) \dxdt  =  \left[  \int_{\Om} \vr \vu \cdot \vc{\phi} \dx  \right]_{t=0}^{t=\tau}    \label{pass-13}
\end{align}
for any $\tau \in [0,T]$ and any $\vc{\phi} \in C_c^{1}( [0,T] \times \Om;\R^2)$. Obviously, the central task left in (\ref{pass-13}) is to verify
\beq \label{pass-14}
\overline{\vr^{\gamma}} =\vr^{\gamma},\,\,\,\, \overline{b^2}=b^2\,\,\,\, \text{ a.e. in } (0,T)\times \Om,
\eeq
which is equivalent to the pointwise convergence of $\{\vr_{\delta}\}_{\delta>0}$ and $\{b_{\delta}\}_{\delta>0}$. This will be the aim of the next subsection. Notice that (\ref{pass-8}) allows to pass to the limits in the equation of continuity (\ref{pass-3}) and magnetic field (\ref{pass-5}) to conclude
\begin{align}
    &  \int_0^{\tau} \int_{\Om} \Big(\vr \p_t \phi+ \vr \vu \cdot \Grad \phi  \Big) \dxdt-\int_0^{\tau} \int_{ \Gamma^{ \rm{in} } }   \vr_{B}  \vu_{B}\cdot \mathbf{n} \phi \dS \dt    \nonumber      \\
    &   \quad  = \int_0^{\tau} \int_{ \Gamma^{ \rm{out} } }   \vr \vu_{B}\cdot \mathbf{n} \phi \dS \dt + \left[  \int_{\Om} \vr \phi \dx  \right]_{t=0}^{t=\tau} , \label{pass-15}  \\
     &  \int_0^{\tau} \int_{\Om} \Big(b \p_t \phi+ b \vu  \cdot \Grad \phi \Big) \dxdt-\int_0^{\tau} \int_{ \Gamma^{ \rm{in} } }   b_{B}  \vu_{B}\cdot \mathbf{n} \phi \dS \dt   \nonumber       \\
    &   \quad  = \int_0^{\tau} \int_{ \Gamma^{ \rm{out} } }   b \vu_{B}\cdot \mathbf{n} \phi \dS \dt + \left[  \int_{\Om} b \phi \dx  \right]_{t=0}^{t=\tau}  \label{pass-16}
\end{align}
for any $\tau \in [0,T]$ and any $\phi\in C^{1}( [0,T] \times \overline{\Om})$.

Before ending this subsection, we report the following useful lemma from Kwon and Novotn\'{y} \cite{KwNo21}.
\begin{Lemma}\label{boundary-esti}
Let $\Om \subset \R^2$ be a bounded Lipschitz domain and $\hat{U}^{-}_{h}$ be an inner neighborhood of its boundary defined as $\hat{U}^{-}_{h}:= \{  x\in \Om\,|\, \text{dist}(x,\p \Om)<h      \}$.
Consider a sequence $\{p_{\delta},\mathbf{z}_{\delta},\mathbf{F}_{\delta},\mathbb{G}_{\delta}  \}_{\delta>0}$ of functions subject to
\[
\p_t \mathbf{z}_{\delta}+ \mathbf{F}_{\delta}+ \Div \mathbb{G}_{\delta} +\Grad p_{\delta}=0\,\,  \text{   in  } \mathcal{D}'((0,T)\times \Om;\R^2).
\]
For $\alpha,\kappa>1$, assume further that
\begin{align*}
    &    \{\mathbf{z}_{\delta}\}_{\delta>0}  \text{  is bounded in  } L^{\infty}(0,T;L^{\alpha}(\Om;\R^2)),        \\
    &   \{\mathbf{F}_{\delta}\}_{\delta>0}  \text{  is bounded in  } L^{\kappa}( (0,T)\times \Om;\R^2      ),   \\
    &   \{\mathbb{G}_{\delta}\}_{\delta>0}  \text{  is bounded in  } L^{\kappa}( (0,T)\times \Om;\R^{2\times 2}      ),  \\
    &  \{p_{\delta}\}_{\delta>0}  \text{  is bounded in  } L^{1}((0,T)\times \Om),
\end{align*}
uniformly with respect to $\delta$.

Then there exists $h_0>0$ and $c=c(\kappa,T,\Om)>0$ such that for any $0<h<h_0$
\[
\int_0^T \int_{ \hat{U}^{-}_{h} } p_{\delta} \dxdt \leq c h^{\Gamma}, \,\,\,\, \Gamma:=\min\{ 1/\alpha',1/\kappa' \}
\]
uniformly with respect to $\delta$.
\end{Lemma}

A detailed proof of this lemma can be seen in Kwon and Novotn\'{y} \cite{KwNo21}. Lemma \ref{boundary-esti} plays a crucial role in passing to the limit in the approximate energy inequality (\ref{pass-6}).

\subsection{Strong convergence of $\{\vr_{\delta}\}_{\delta>0}$ and $\{b_{\delta}\}_{\delta>0}$ }
As observed in the previous subsection, the key issue lies in showing the strong convergence of $\{\vr_{\delta}\}_{\delta>0}$ and $\{b_{\delta}\}_{\delta>0}$ due to the nonlinearies from the pressure and magnetic field. To this end, we first set
\[
\zeta_{\delta}:=\left\{
\begin{aligned}
&  \f{b_{\delta}}{\vr_{\delta}}  \text{  if  } \vr_{\delta}>0, \\
&  \f{C_{\ast}+C^{\ast}}{2}  \text{  if  } \vr_{\delta}=0 ,
\end{aligned}
\right. \,\,\,\, \,\,\,\,
\zeta:=\left\{
\begin{aligned}
&  \f{b}{\vr}  \text{  if  } \vr>0, \\
&  \f{C_{\ast}+C^{\ast}}{2}  \text{  if  } \vr=0 .
\end{aligned}
\right.
\]
It turns out that $\{\zeta_{\delta} \}_{\delta>0}$ enjoys the \emph{almost compactness} property; it dates back to the renormalization and transport theory of DiPerna and Lions \cite{DL89}. The theory was later adapted to
the compressible Navier-Stokes system \cite{FNP} and more recently in compressible two-fluid models \cite{VWY,NP} with homogeneous Dirichlet boundary conditions. Very recently, Kra\v{c}mar et al. \cite{KKNN22} extended the theory to the families of transport equations with general boundary conditions. In our context of compressible MHD model, it holds that
\begin{Lemma}\label{almost-comp}
\begin{align*}
    &  \lim_{\delta \rightarrow 0}  \int_{\Om} \vr_{\delta} |\zeta_{\delta}-\zeta|^2 (\tau)\dx=0,          \\
    &   \lim_{\delta \rightarrow 0} \int_0^{\tau} \int_{  \Gamma^{ \rm{out} } }  \vr_{\delta} |\zeta_{\delta}-\zeta|^2   \vu_{B} \cdot \mathbf{n}   \dS\dt =0
\end{align*}
for any $\tau \in [0,T]$.
\end{Lemma}
{\emph {Proof of Lemma \ref{almost-comp}.}} Essentially based on the renormalization technique for transport equations with general boundary data, see Corollary 5 in Kra\v{c}mar et al. \cite{KKNN22}, it follows from the uniform estimates (\ref{pass-7}), the approximate equation of continuity (\ref{pass-3}) and magnetic field (\ref{pass-5}) that $(\vr_{\delta} \zeta_{\delta}^2,\vu_{\delta})$ satisfies the equation of continuity\footnote{Due to the nice integrability from the magnetic field and the domination condition, the density naturally inherits a good integrability; this perfectly matches the hypotheses when applying the renormalization theory from Kra\v{c}mar et al. \cite{KKNN22} and allows to prove existence of weak solutions for any $\gamma>1$. For more discussions on the domination condition and adiabatic exponent, see Section \ref{con-rem}.}, i.e.,
\begin{align}
    &   \int_{\Om}  \vr_{\delta} \zeta_{\delta}^2 (\tau) \dx + \int_0^{\tau} \int_{  \Gamma^{ \rm{out} } }  \vr_{\delta} \zeta_{\delta}^2  \vu_{B} \cdot \mathbf{n}   \dS\dt     \nonumber        \\
    &   \quad  =\int_{\Om}  \vr_{\delta} \zeta_{\delta}^2 (0) \dx  - \int_0^{\tau} \int_{  \Gamma^{ \rm{in} } }  \vr_{B} \zeta_{B}^2  \vu_{B} \cdot \mathbf{n}   \dS\dt  , \,\, \,\, \zeta_{B}:= b_{B}/\vr_{B}   .  \label{pass-17}
\end{align}
In the same spirit, it follows from (\ref{pass-15})-(\ref{pass-16}) that
\begin{align}
    &   \int_{\Om}  \vr \zeta^2 (\tau) \dx + \int_0^{\tau} \int_{  \Gamma^{ \rm{out} } }  \vr \zeta^2  \vu_{B} \cdot \mathbf{n}   \dS\dt      \nonumber         \\
    &   \quad  =\int_{\Om}  \vr \zeta^2 (0) \dx  - \int_0^{\tau} \int_{  \Gamma^{ \rm{in} } }  \vr_{B} \zeta_{B}^2  \vu_{B} \cdot \mathbf{n}   \dS\dt  , \,\, \,\, \zeta_{B}:= b_{B}/\vr_{B}   .  \label{pass-18}
\end{align}

Observe that
\begin{align}
    &    \lim_{\delta \rightarrow 0}  \left(   \int_{\Om} \vr_{\delta} |\zeta_{\delta}-\zeta|^2 (\tau)\dx +     \int_0^{\tau} \int_{  \Gamma^{ \rm{out} } }  \vr_{\delta} |\zeta_{\delta}-\zeta|^2   \vu_{B} \cdot \mathbf{n}   \dS\dt        \right)         \nonumber    \\
    &   \quad  =  \lim_{\delta \rightarrow 0}   \left(   \int_{\Om} \vr_{\delta} \zeta_{\delta}^2 (\tau)\dx  +  \int_0^{\tau} \int_{  \Gamma^{ \rm{out} } }     \vr_{\delta} \zeta_{\delta}^2  \vu_{B} \cdot \mathbf{n}   \dS\dt       \right)          \nonumber     \\
    &   \quad  \quad   \quad    - \left(    \int_{\Om} \vr \zeta^2 (\tau) \dx +    \int_0^{\tau} \int_{  \Gamma^{ \rm{out} } } \vr \zeta^2     \vu_{B} \cdot \mathbf{n}   \dS\dt        \right) , \label{pass-19}
\end{align}
where we used the basic fact that
\begin{align*}
    &   \lim_{\delta \rightarrow 0}    \int_{\Om}  \vr_{\delta} \zeta_{\delta} \zeta(\tau) \dx= \lim_{\delta \rightarrow 0}    \int_{\Om} b_{\delta} \zeta(\tau) \dx =  \int_{\Om} b \zeta(\tau) \dx
    =   \int_{\Om} \vr \zeta^2(\tau) \dx,      \\
    &   \lim_{\delta \rightarrow 0}  \int_0^{\tau}  \int_{\Om}  \vr_{\delta} \zeta_{\delta} \zeta \vu_{B} \cdot \mathbf{n}   \dS\dt  = \lim_{\delta \rightarrow 0} \int_0^{\tau}   \int_{\Om} b_{\delta} \zeta \vu_{B} \cdot \mathbf{n}   \dS\dt  =\int_0^{\tau}
     \int_{\Om} b \zeta \vu_{B} \cdot \mathbf{n}   \dS\dt    \\
    &   \quad     =\int_0^{\tau}
     \int_{\Om} \vr \zeta^2 \vu_{B} \cdot \mathbf{n}   \dS\dt .
\end{align*}
Combining (\ref{pass-17})-(\ref{pass-19}), we arrive at Lemma \ref{almost-comp}.          $\Box$

As a corollary of Lemma \ref{almost-comp}, we easily infer, through a routine interpolation argument, from the uniform estimates (\ref{pass-7}) and the improved integrability estimate (\ref{pass-11}) that
\begin{Corollary}\label{Cor4.1}
\begin{align*}
    &  \lim_{\delta \rightarrow 0}  \int_{\Om} \vr_{\delta}^{p} |\zeta_{\delta}-\zeta|^{q} (\tau)\dx=0,  \,\, \text{  for any  }  1\leq p < \max\{ 2,\gamma \} ,1\leq q <\infty,     \\
    &   \lim_{\delta \rightarrow 0} \int_0^{T} \int_{  K}  \vr_{\delta}^{p}  |\zeta_{\delta}-\zeta|^{q}     \dx\dt =0,  \,\, \text{  for any  }  1\leq p < \max\{ \gamma+\vt,2+\vt \}   ,1\leq q <\infty,
\end{align*}
for any compact $K\subset \Om$. Here, $\vt>0$ depends solely on $\gamma$ and comes from (\ref{pass-11}).
\end{Corollary}

To proceed, we make use of the following decomposition
\begin{align}
    \vr_{\delta}^{\gamma}+\f{1}{2}b_{\delta}^2 &  = \vr_{\delta}^{\gamma}+ \f{1}{2}  \zeta_{\delta}^2  \vr_{\delta}^2      \nonumber     \\
    &   = \left( \vr_{\delta}^{\gamma}+ \f{1}{2}  \zeta_{\delta}^2  \vr_{\delta}^2 \right)-  \left( \vr_{\delta}^{\gamma}+ \f{1}{2}  \zeta^2  \vr_{\delta}^2 \right)+ \left( \vr_{\delta}^{\gamma}+ \f{1}{2}  \zeta^2  \vr_{\delta}^2 \right)    \nonumber    \\
    &  =  \f{1}{2}  \vr_{\delta}^2 ( \zeta_{\delta}^2 - \zeta^2  )    +\left( \vr_{\delta}^{\gamma}+ \f{1}{2}  \zeta^2  \vr_{\delta}^2 \right)  ;   \label{pass-20}
\end{align}
whence with the help of Corollary \ref{Cor4.1} and the available convergence results, we may pass to the limit in the approximate momentum equation (\ref{pass-4}) to conclude a refined version of (\ref{pass-13})
\begin{align}
    &  \int_0^{\tau} \int_{\Om} \Big(   \vr \vu \cdot \p_t \vc{\phi} +\vr \vu \otimes \vu :\Grad \vc{\phi}+ \overline{  \overline{   \mathcal{P}(\vr,\zeta)      }    } \,  \Div \vc{\phi}     \nonumber    \\
    &   \quad \quad     -\mathbb{S}(\Grad \vu ):\Grad \vc{\phi}                     \Big) \dxdt  =  \left[  \int_{\Om} \vr \vu \cdot \vc{\phi} \dx  \right]_{t=0}^{t=\tau}    \label{pass-21}
\end{align}
for any $\tau \in [0,T]$ and any $\vc{\phi} \in C_c^{1}( [0,T] \times \Om;\R^2)$. Here and in the sequel, we denote by $\overline{  \overline{   \mathcal{P}(\vr,\zeta)      }    }$ the weak limit of
$\{  \vr_{\delta}^{\gamma}+ \f{1}{2}  \zeta^2  \vr_{\delta}^2 \}_{\delta>0}$.

We are now in a position to employ the remarkable effective viscous flux identity. To do this, we introduce a sequence of cut-off functions as in \cite{FNP}
\[
T_k(z):= k    T\left( \f{z}{k}  \right),\,\,k \geq 1,
\]
where $T(z)$ is a smooth concave function in $[0,\infty)$ subject to
\[
T_k(z)= \left\{
\begin{aligned}
&  z  \text{  if  } z\in [0,1], \\
&  2  \text{  if  } z\in [3,\infty) .
\end{aligned}
\right.
\]

\begin{Lemma}\label{effec-flux}
It holds that
\begin{align*}
    &   \lim_{\delta \rightarrow 0} \int_0^T \eta \int_{\Om} \phi \left\{    \left( \vr_{\delta}^{\gamma}+ \f{1}{2}  \zeta^2  \vr_{\delta}^2 \right)-(\lambda+2\mu) \Div \vu_{\delta}                \right\} T_k(\vr_{\delta})  \dxdt          \\
    &   \quad  =   \int_0^T \eta \int_{\Om} \phi    \left\{       \overline{  \overline{   \mathcal{P}(\vr,\zeta)      }    }       -(\lambda+2\mu)  \Div \vu    \right\}  \overline{T_k(\vr)}   \dxdt
\end{align*}
for any $\eta \in C_c^{\infty}( (0,T) ),\phi \in C_c^{\infty}( \Om )$.
\end{Lemma}

We remark that the proof of Lemma \ref{effec-flux} shares similar features from the compressible Navier-Stokes system \cite{FNP,LS}, the compressible two-fluid models \cite{NP,VWY,KKNN22} and the compressible MHD models \cite{LS19,LS21}. Specifically, one first uses $\eta \phi \mathcal{A}[1_{\Om} T_k(\vr_{\delta}) ]$ as a test function in the approximate momentum equation (\ref{pass-4}), where $\mathcal{A}=(\mathcal{A}_1,\mathcal{A}_2)$ with $\mathcal{A}_j=\Delta^{-1}\p_{x_j} $ for $j=1,2$.
Next, we test the limit momentum equation (\ref{pass-21}) via $\eta \phi \mathcal{A}[1_{\Om} T_k(\vr) ]$. It remains to pass to the limits in the resulting approximate momentum equation (\ref{pass-4}) with test function $\eta \phi \mathcal{A}[1_{\Om} T_k(\vr_{\delta}) ]$ and combine with the limit momentum equation (\ref{pass-21}) with test function $\eta \phi \mathcal{A}[1_{\Om} T_k(\vr) ]$, employing conveniently in particular the renormalization technique for continuity equation, the H\"{o}rmander-Mikhlin theorem and the Div-Curl lemma. We do not intend to give the details here since they can be carried out in a routine manner, see for instance Lemma 4.1 in \cite{LS19} or Proposition 14 in \cite{KKNN22}. Notice that the decomposition of pressure term (\ref{pass-20}) and Corollary \ref{Cor4.1} play crucial role in handling the term
\[
\lim_{\delta \rightarrow 0}   \int_0^T \eta \int_{\Om}  \phi  \left( \vr_{\delta}^{\gamma}+\f{1}{2}b_{\delta}^2 \right) T_k(\vr_{\delta})   \dxdt          .
\]

With Lemma \ref{effec-flux} at hand, the remaining steps follow in a similar way as the compressible Navier-Stokes system, based on a careful analysis of oscillations defect measure and the renormalization technique.
These are nowadays standard in the theory of compressible fluids and we refer to sections 7.4-7.5 in \cite{KKNN22}, sections 4.4-4.6 in \cite{FNP} for the details. Consequently, we arrive at
\[
\int_{\Om}  ( \overline{\vr \log \vr}-   \vr \log \vr    )(\tau) \dx  +   \int_0^{\tau} \int_{  \Gamma^{ \rm{out} } }   ( \overline{\vr \log \vr}-   \vr \log \vr    ) \vu_{B} \cdot \mathbf{n}   \dS\dt
\leq 0;
\]
whence the strict convexity of $z \mapsto z\log z$ in $[0,\infty)$ yields
\beq \label{pass-22}
\vr_{\delta} \rightarrow \vr \text{  a.e. in  } (0,T)\times \Om , \,\, \text{  a.e. in  } (0,T) \times  \Gamma^{ \rm{out} }.
\eeq
It follows readily that
\beq \label{pass-23}
b_{\delta} \rightarrow b \text{  a.e. in  } (0,T)\times \Om , \,\, \text{  a.e. in  } (0,T) \times  \Gamma^{ \rm{out} }.
\eeq
upon seeing the basic relation
\begin{align*}
    |b_{\delta} -b |   & = | \vr_{\delta}\zeta_{\delta}-\vr \zeta  |= |      \vr_{\delta}\zeta_{\delta}   -\vr_{\delta}\zeta+\vr_{\delta}\zeta  -   \vr \zeta         |           \\
    &   \leq  \vr_{\delta}    |\zeta_{\delta}  -\zeta|+ \zeta |\vr_{\delta}-\vr|
\end{align*}
together with Corollary \ref{Cor4.1} and (\ref{pass-22}). Thus,
\[
 \overline{  \overline{   \mathcal{P}(\vr,\zeta)      }    } = \vr^{\gamma}+\f{1}{2} b^2
\]
and the weak formulation (\ref{pass-21}) then reduces to the desired balance of momentum in sense of distributions.

It remains to pass to the limit in the approximate energy inequality (\ref{pass-6}). Observe that all the terms could be easily handled via the available convergences with the only exception
\[
 \int_0^{\tau} \int_{\Om}  \Big(  \vr_{\delta}^{\gamma}+\f{1}{2}b_{\delta}^2+\delta(\vr_{\delta}+b_{\delta})^{\beta}   \Big)           \Div \vu_B \dxdt ,
\]
since the improved integrability estimate (\ref{pass-10}) holds merely on compact subsets of $\Om$. To this end, for sufficiently small $h>0$ we make the decomposition
\begin{align}
    &     \int_0^{\tau} \int_{\Om}  \Big(  \vr_{\delta}^{\gamma}+\f{1}{2}b_{\delta}^2+\delta(\vr_{\delta}+b_{\delta})^{\beta}   \Big)           \Div \vu_B \dxdt    \nonumber    \\
    &   \quad  =    \int_0^{\tau} \int_{ \Om \backslash \hat{U}^{-}_{h}   }  \Big(  \vr_{\delta}^{\gamma}+\f{1}{2}b_{\delta}^2+\delta(\vr_{\delta}+b_{\delta})^{\beta}   \Big)           \Div \vu_B \dxdt  \nonumber    \\
    &   \quad  \quad  \quad +   \int_0^{\tau} \int_{\hat{U}^{-}_{h}   }  \Big(  \vr_{\delta}^{\gamma}+\f{1}{2}b_{\delta}^2+\delta(\vr_{\delta}+b_{\delta})^{\beta}   \Big)           \Div \vu_B \dxdt .   \label{pass-24}
\end{align}
For fixed $h>0$, the first integral on the right-hand side of (\ref{pass-24}) tends to
\[
 \int_0^{\tau} \int_{ \Om \backslash \hat{U}^{-}_{h}   }  \Big(  \vr^{\gamma}+\f{1}{2}b^2  \Big)           \Div \vu_B \dxdt
\]
as $\delta \rightarrow 0$ mainly due to (\ref{pass-12}) and the strong convergences (\ref{pass-22})-(\ref{pass-23}). Thanks to Lemma \ref{boundary-esti}, the second integral on the right-hand side of (\ref{pass-24}) is bounded
uniformly by $c h^{\Gamma}$ for fixed $h>0$. Upon passing to the limit $h\rightarrow 0$, we obtain the desired integral.

Up to now, we have proved the existence of global weak solutions with initial-boundary data satisfying the hypotheses of Lemma \ref{lem-3.1}. Observe that the general case of finite energy initial data in Theorem \ref{TH1} is treated exactly in the same manner, with a suitable regularization of initial data parameterized by $\delta$; see for instance section 4 of \cite{FNP} or section 3 of \cite{LS19}. The proof of Theorem \ref{TH1} is completely finished.

\section{Weak-strong uniqueness principle}\label{weak-strong}
This section is devoted to the stability of strong solution in the class of weak solutions. Our strategy is based on the well-known relative entropy method, see \cite{FJA12,FNS11,KwNo21} in the context of compressible Navier-Stokes system.

\subsection{Relative energy inequality}\label{relative-energy}

Let $(\vr,\vu,b)$ and $(\tilde{\vr},\tilde{\vu},\tilde{b}  )$ be as in Theorem \ref{TH2}. For simplicity, we set $\tilde{\vv}:=\tilde{\vu}-\vu_{B}$.
To this end, we define the relative energy as follows
\begin{align*}
    &   \mathcal{E}(\vr,\vu,b\,|\, \tilde{\vr},\tilde{\vu},\tilde{b}  ):= \int_{\Om} \left(  \f{1}{2} \vr |\vv-\tilde{\vv}|^2 +H(\vr)-H(\tilde{\vr})-H'(\tilde{\vr})(\vr-\tilde{\vr})+\f{1}{2}|b-\tilde{b}|^2 \right) \dx  ,
\end{align*}
where we have set $\tilde{\vv}:=\tilde{\vu}-\vu_{B} $ and the pressure potential $H(\vr):=\f{ \vr^{\gamma} }{\gamma-1}$ for simplicity. Observe that we may reformulate the relative energy as
\begin{align*}
   \mathcal{E}(\vr,\vu,b\,|\, \tilde{\vr},\tilde{\vu},\tilde{b}  )  &  = \int_{\Om}  \left(  \f{1}{2} \vr |\vv|^2+H(\vr)+\f{1}{2}  b^2    \right) \dx - \int_{\Om}  \vr \vv \cdot  \tilde{\vv} \dx              \\
    &   \quad  +  \int_{\Om}  \f{1}{2} \vr | \tilde{\vv} |^2  \dx - \int_{\Om}  H'(\tilde{\vr})  \vr  \dx   -     \int_{\Om}    b  \tilde{b}   \dx    \\
    &   \quad  + \int_{\Om}  \left(  p(\tilde{\vr}) +  \f{1}{2} |\tilde{b} |^2      \right)         \dx     \\
    & =: \sum_{j=1}^6 I_j.
\end{align*}
It remains to estimate $I_j,j=1,...,6$ suitably. For $I_1$, we directly invoke the energy inequality (\ref{energy-inequa}) to obtain
\begin{align}
    &    I_1   +  \int_0^{\tau} \int_{\Om}\mathbb{S}(\Grad \vu):\Grad \vu \dxdt   +  \int_0^{\tau} \int_{ \Gamma^{ \rm{out} } }  \left(  H(\vr) +\f{1}{2} b^2  \right) \vu_{B} \cdot \mathbf{n}  \dS \dt   \nonumber     \\
    &   \quad     \leq      \int_{\Om}    \left( \f{1}{2} \vr_0 | \vv_0|^2 +H(\vr_0) +\f{1}{2} b_0^2    \right)  \dx    -\int_0^{\tau} \int_{ \Gamma^{ \rm{in} } }  \left(  H(\vr_{B}) +\f{1}{2} b_{B}^2  \right) \vu_{B} \cdot \mathbf{n}  \dS \dt   \label{weak-strong-1}           \\
    & \quad  \quad + \int_0^{\tau} \int_{\Om}  \left[   -\left( p(\vr)+\f{1}{2}b^2  \right) \Div \vu_{B} -\vr \vu \cdot \Grad \vu_{B} \cdot (\vu-\vu_{B}) +  \mathbb{S}(\Grad \vu):\Grad \vu_{B}                    \right] \dxdt .   \nonumber
\end{align}
For $I_2$, we test the momentum equation (\ref{balan-momentu}) by $-\tilde{\vv}$ to obtain
\begin{align}
   I_2 +  \int_{\Om} \vr_0 \vv_0 \cdot \tilde{\vv}(0) \dx          & =   \int_0^{\tau} \int_{\Om}  \Big( -\vr \vv \cdot \p_t \tilde{\vv} -\vr \vu \cdot \Grad \tilde{\vv} \cdot \vu+\vr \vu \cdot \Grad (\vu_{B} \cdot \tilde{\vv})               \Big) \dxdt          \nonumber     \\
   &    \quad -  \int_0^{\tau} \int_{\Om} \left( p(\vr)+\f{1}{2}b^2 \right)\Div  \tilde{\vv} \dxdt+    \int_0^{\tau}  \int_{\Om}  \mathbb{S}(\Grad \vu):\Grad    \tilde{\vv}        \dxdt  .        \label{weak-strong-2}
\end{align}
Using $\f{1}{2} |\tilde{\vv}|^2$ as a test function in the equation of continuity (\ref{equ-conti}),
\beq\label{weak-strong-3}
I_3 -  \int_{\Om} \f{1}{2} \vr_0 |\tilde{\vv}(0) |^2 \dx = \int_0^{\tau} \int_{\Om} \Big(  \vr \tilde{\vv} \cdot \p_t \tilde{\vv} + \vr \vu \cdot \Grad \tilde{\vv}  \cdot \tilde{\vv}      \Big) \dxdt.
\eeq
Next, we choose $-H'(\tilde{\vr})$ as a test function in the equation of continuity (\ref{equ-conti}) and notice that $H''(\tilde{\vr})=\f{p'(\tilde{\vr})}{\tilde{\vr}}$ in order to see
\begin{align}
   I_4 +  \int_{\Om} \vr_0   H'(  \tilde{\vr}(0) )   \dx          & = -  \int_0^{\tau} \int_{\Om}   \f{p'(\tilde{\vr})}{\tilde{\vr}  }  \Big(  \vr \p_t \tilde{\vr} +\vr \vu \cdot \Grad \tilde{\vr}       \Big)                 \dxdt          \nonumber     \\
   &    \quad +  \int_0^{\tau} \int_{ \Gamma^{ \rm{in} } }   \vr_{B}  \vu_{B}\cdot \mathbf{n}  H'(\tilde{\vr})  \dS \dt        \label{weak-strong-4} \\
   &      \quad +        \int_0^{\tau} \int_{ \Gamma^{ \rm{out} } }   \vr  \vu_{B}\cdot \mathbf{n}  H'(\tilde{\vr})  \dS \dt  .         \nonumber
\end{align}
For $I_5$, we test the equation of magnetic field (\ref{equ-magnetic}) by $-\tilde{b}$ to arrive at
\begin{align}
   I_5 +  \int_{\Om} b_0     \tilde{b}(0)    \dx          & = -  \int_0^{\tau} \int_{\Om}     \Big(  b \p_t \tilde{b} + b \vu \cdot \Grad \tilde{b}       \Big)                 \dxdt          \nonumber     \\
   &    \quad +  \int_0^{\tau} \int_{ \Gamma^{ \rm{in} } }   b_{B}^2  \vu_{B}\cdot \mathbf{n}    \dS \dt        \label{weak-strong-5} \\
   &      \quad +        \int_0^{\tau} \int_{ \Gamma^{ \rm{out} } }   b \vu_{B}\cdot \mathbf{n}  \tilde{b}  \dS \dt  .         \nonumber
\end{align}
To estimate $I_6$, we calculate
\begin{align}
   I_6          & =   \int_{\Om}  p(  \tilde{\vr} )  \dx - \int_{\Om}  p(  \tilde{\vr}(0) )  \dx   + \int_{\Om}  p(  \tilde{\vr}(0) )  \dx                   \nonumber     \\
   &    \quad +    \int_{\Om} \f{1}{2}  |\tilde{b} |^2 \dx  -\int_{\Om} \f{1}{2}  |\tilde{b}(0) |^2 \dx +\int_{\Om} \f{1}{2}  |\tilde{b}(0) |^2 \dx            \nonumber      \\
   &     =    \int_0^{\tau} \int_{\Om} p'(  \tilde{\vr} )   \p_t \tilde{\vr}  \dxdt + \int_{\Om}  \Big(   H'( \tilde{\vr}(0) )\tilde{\vr}(0)-H( \tilde{\vr}(0) )    \Big)             \dx                                            \nonumber \\
   &        \quad +  \int_0^{\tau} \int_{\Om}      \tilde{b}  \p_t  \tilde{b}  \dxdt+   \int_{\Om} \f{1}{2}  |\tilde{b}(0) |^2 \dx.                 \label{weak-strong-6}
\end{align}

Summing up (\ref{weak-strong-1})-(\ref{weak-strong-6}) and using the fact that $(\vr,\vu,b)$ and $(\tilde{\vr},\tilde{\vu},\tilde{b}  )$ share the same initial-boundary data, one obtains
\begin{align}
   &      \mathcal{E}(\vr,\vu,b\,|\, \tilde{\vr},\tilde{\vu},\tilde{b}  )+  \int_0^{\tau} \int_{ \Gamma^{ \rm{out} } }  \Big(    H(\vr)-H'( \tilde{\vr})\vr               \Big)  \vu_{B}\cdot \mathbf{n}  \dS \dt          \nonumber     \\
 &      \quad\quad +    \int_0^{\tau} \int_{ \Gamma^{ \rm{out} } }  \Big(      \f{1}{2}b^2 -b \tilde{b}       \Big)  \vu_{B}\cdot \mathbf{n}  \dS \dt           \nonumber     \\
 & \quad\quad +      \int_0^{\tau} \int_{\Om}          \mathbb{S}(\Grad \vu):\Grad (\vv-\tilde{\vv}) \dxdt                     \nonumber     \\
 &   \quad \leq   \int_0^{\tau} \int_{ \Gamma^{ \rm{in} } }   \Big(    -H(\vr_{B})+ H'(\tilde{\vr}) \vr_{B}     \Big)  \vu_{B}\cdot \mathbf{n}  \dS \dt + \int_0^{\tau} \int_{ \Gamma^{ \rm{in} } } \f{1}{2}   b_{B}^2  \vu_{B}\cdot \mathbf{n}   \dS \dt             \nonumber     \\
 & \quad\quad +     \int_0^{\tau} \int_{\Om}  \Big(  \vr (\tilde{\vv}-\vv)\p_t \tilde{\vv} +\vr\vu \cdot \Grad \tilde{\vu} \cdot (\tilde{\vv} -\vv)      \Big)      \dxdt              \label{weak-strong-7}      \\
  & \quad\quad +    \int_0^{\tau} \int_{\Om} \Big(  p(\tilde{\vr})-p(\vr)        \Big)  \Div \tilde{\vu}  \dxdt   +    \int_0^{\tau} \int_{\Om} \Big(  \f{1}{2} |\tilde{b} |^2 -\f{1}{2} b^2      \Big)  \Div \tilde{\vu}  \dxdt                       \nonumber     \\
   & \quad\quad +   \int_0^{\tau} \int_{\Om} \left(            \f{\tilde{\vr}-\vr   }{  \tilde{\vr}  }p'(\tilde{\vr}) (\p_t \tilde{\vr} +\vu \cdot \Grad \tilde{\vr})-p'(\tilde{\vr}) \vu \cdot \Grad \tilde{\vr}-p(\tilde{\vr})
   \Div \tilde{\vu}          \right) \dxdt                                                  \nonumber     \\
   &  \quad \quad    -  \int_0^{\tau} \int_{\Om} \Big(  (\tilde{b}-b)(\p_t \tilde{b} +\vu \cdot \Grad \tilde{b})   -\tilde{b} \vu \cdot \Grad \tilde{b} -\f{1}{2} |\tilde{b}|^2 \Div \tilde{\vu}     \Big)         \dxdt .   \nonumber
\end{align}
Next we shall simplify the expressions on the right-hand side of (\ref{weak-strong-7}) suitably. Based on the fact that $\p_t \tilde{\vr} + \Div(\tilde{\vr}\tilde{\vu})=0 $, a direct calculation gives
\begin{align}
   &        \int_0^{\tau} \int_{\Om} \left(            \f{\tilde{\vr}-\vr   }{  \tilde{\vr}  }p'(\tilde{\vr}) (\p_t \tilde{\vr} +\vu \cdot \Grad \tilde{\vr})-p'(\tilde{\vr}) \vu \cdot \Grad \tilde{\vr}-p(\tilde{\vr})
   \Div \tilde{\vu}          \right) \dxdt              \nonumber     \\
    &   \quad =   \int_0^{\tau} \int_{\Om} \left(    \f{\tilde{\vr}-\vr   }{  \tilde{\vr}  }p'(\tilde{\vr})   (\vu-\tilde{\vu}) \cdot \Grad \tilde{\vr}   +  (\vr-\tilde{\vr}) p'(\tilde{\vr}) \Div \tilde{\vu}
    -p'(\tilde{\vr}) \vv \cdot \Grad \tilde{\vr} -p(\tilde{\vr}) \Div \tilde{\vv}   \right)  \dxdt       \nonumber           \\
    &   \quad\quad -   \int_0^{\tau} \int_{ \Gamma^{ \rm{in} } }  p(\vr_{B})      \vu_{B}\cdot \mathbf{n}  \dS \dt - \int_0^{\tau} \int_{ \Gamma^{ \rm{out} } }  p(\tilde{\vr})      \vu_{B}\cdot \mathbf{n}  \dS \dt.                                 \label{weak-strong-8}
    \end{align}
In the same spirit, we make use of the magnetic equation satisfied by $(\tilde{b},\tilde{\vu})$ to infer
\begin{align}
&   \int_0^{\tau} \int_{\Om} \Big(  (\tilde{b}-b)(\p_t \tilde{b} +\vu \cdot \Grad \tilde{b})   -\tilde{b} \vu \cdot \Grad \tilde{b} -\f{1}{2} |\tilde{b}|^2 \Div \tilde{\vu}     \Big)         \dxdt             \nonumber     \\
&    \quad =  \int_0^{\tau} \int_{\Om}   \Big(   (\tilde{b}-b)(   \vu- \tilde{\vu}  ) \cdot \Grad \tilde{b}+ (b- \tilde{b}) \tilde{b} \Div \tilde{\vu}-\tilde{b} \vv \cdot \Grad \tilde{b} -\f{1}{2}|\tilde{b}|^2 \Div \tilde{\vv}                    \Big)   \dxdt                             \nonumber     \\
&     \quad\quad -   \int_0^{\tau} \int_{ \Gamma^{ \rm{in} } }      \f{1}{2} b_{B}^ 2     \vu_{B}\cdot \mathbf{n}  \dS \dt -  \int_0^{\tau} \int_{ \Gamma^{ \rm{out} } } \f{1}{2} |\tilde{b} |^2    \vu_{B}\cdot \mathbf{n}   \dS \dt.   \label{weak-strong-9}
\end{align}
Multiplying the momentum equation
\[
\tilde{\vr} (\p_t \tilde{\vu} +\tilde{\vu} \cdot \Grad \tilde{\vu} ) + \Grad \left(p(\tilde{\vr})+\f{1}{2}|\tilde{b}|^2         \right)= \Div \mathbb{S}(\Grad \tilde{\vu} )
\]
by $\vv-\tilde{\vv}$ and integration by parts give rise to
\begin{align}
- \int_0^{\tau} \int_{\Om}   \mathbb{S}(\Grad \tilde{\vu} ) : \Grad (\vv- \tilde{\vv}) \dxdt & =   \int_0^{\tau} \int_{\Om}  \tilde{\vr} (\p_t \tilde{\vu} +\tilde{\vu} \cdot \Grad \tilde{\vu}  ) \cdot (\vv-\tilde{\vv})  \dxdt                   \nonumber     \\
&    \quad +    \int_0^{\tau} \int_{\Om}   p'(\tilde{\vr}) \vv\cdot \Grad \tilde{\vr}           \dxdt     +       \int_0^{\tau} \int_{\Om}     p(\tilde{\vr}) \Div \tilde{\vv}  \dxdt       \label{weak-strong-10} \\
&    \quad +     \int_0^{\tau} \int_{\Om}    \tilde{b} \vv \cdot \Grad \tilde{b} \dxdt  +    \int_0^{\tau} \int_{\Om}   \f{1}{2} | \tilde{b} |^2 \Div  \tilde{\vv}       \dxdt .                     \nonumber
\end{align}
Therefore, we first add $- \int_0^{\tau} \int_{\Om}   \mathbb{S}(\Grad \tilde{\vu} ) : \Grad (\vv- \tilde{\vv}) \dxdt$ to both sides of (\ref{weak-strong-7}) and combine the relations (\ref{weak-strong-8})-(\ref{weak-strong-10}) to conclude
\begin{align}
&   \mathcal{E}(\vr,\vu,b\,|\, \tilde{\vr},\tilde{\vu},\tilde{b}  ) +    \int_0^{\tau} \int_{\Om}          \mathbb{S}(\Grad (\vv-\tilde{\vv})):\Grad (\vv-\tilde{\vv}) \dxdt           \nonumber     \\
&  \quad \quad +     \int_0^{\tau} \int_{ \Gamma^{ \rm{out} } }  \Big(    H(\vr)-H(\tilde{\vr})-  H'( \tilde{\vr})(\vr-\tilde{\vr})               \Big)  \vu_{B}\cdot \mathbf{n}  \dS \dt      \nonumber     \\
&  \quad \quad +     \int_0^{\tau} \int_{ \Gamma^{ \rm{out} } }  \Big(    \f{1}{2}b^2  -\f{1}{2}|\tilde{b}|^2 - \tilde{b} (b-\tilde{b})             \Big)  \vu_{B}\cdot \mathbf{n}  \dS \dt      \nonumber     \\
&   \quad  \leq    \int_0^{\tau} \int_{\Om}  \Big(    (\vr-\tilde{\vr} )(\tilde{\vv}-\vv )\p_t \tilde{\vu} +(\vr-\tilde{\vr} )\tilde{\vu} \cdot \Grad \tilde{\vu} \cdot (\tilde{\vv}-\vv )
+ \vr (\vv-\tilde{\vv} ) \cdot \Grad \tilde{\vu}  \cdot (\tilde{\vv}-\vv )   \Big) \dxdt                         \nonumber     \\
&    \quad \quad +  \int_0^{\tau} \int_{\Om}   \left(  1-\f{\vr }{  \tilde{\vr}  }     \right) p'(\tilde{\vr}) (\vv-\tilde{\vv}) \cdot \Grad \tilde{\vr} \dxdt
    +  \int_0^{\tau} \int_{\Om}          (\tilde{b}-b)        (\vv-\tilde{\vv}) \cdot \Grad \tilde{b} \dxdt         \nonumber     \\
    &   \quad \quad   -   \int_0^{\tau} \int_{\Om}    \Big(  p(\vr)-p(\tilde{\vr})  -p'(\tilde{\vr})(\vr-\tilde{\vr})    \Big) \Div \tilde{\vu}  \dxdt          \label{weak-strong-11}    \\
    & \quad \quad    -    \int_0^{\tau} \int_{\Om}    \Big(  \f{1}{2} b^2 -\f{1}{2} |\tilde{b}|^2 -\tilde{b}(b-\tilde{b} )                          \Big) \Div \tilde{\vu}   \dxdt            \nonumber \\
      & \quad \quad     -  \int_0^{\tau} \int_{ \Gamma^{ \rm{in} } }      \Big(   H(\vr_{B})-H(\tilde{\vr})-H'(\tilde{\vr})(\vr_{B}-\tilde{\vr})      \Big)                 \vu_{B}\cdot \mathbf{n}          \dS \dt                             \nonumber \\
      & \quad \quad     -  \int_0^{\tau} \int_{ \Gamma^{ \rm{in} } }      \Big(   \f{1}{2}b_{B}^2  -\f{1}{2}|\tilde{b} |^2      -\tilde{b} (b_{B}-\tilde{b} )         \Big)                 \vu_{B}\cdot \mathbf{n}          \dS \dt.                             \nonumber
\end{align}
Seeing again that $(\vr,b)$ and $(\tilde{\vr},\tilde{b})$ share the same boundary data, the surface integrals on the inflow part from the right-hand side vanish; whence we may furthermore drop the surface integrals on the outflow part from the left-hand side to obtain
\begin{align}
&   \mathcal{E}(\vr,\vu,b\,|\, \tilde{\vr},\tilde{\vu},\tilde{b}  ) +    \int_0^{\tau} \int_{\Om}          \mathbb{S}(\Grad (\vv-\tilde{\vv})):\Grad (\vv-\tilde{\vv}) \dxdt           \nonumber     \\
&   \quad  \leq    \int_0^{\tau} \int_{\Om}  \Big(    (\vr-\tilde{\vr} )(\tilde{\vv}-\vv )\cdot  \p_t \tilde{\vu} +(\vr-\tilde{\vr} )\tilde{\vu} \cdot \Grad \tilde{\vu} \cdot (\tilde{\vv}-\vv )
+ \vr (\vv-\tilde{\vv} ) \cdot \Grad \tilde{\vu}  \cdot (\tilde{\vv}-\vv )   \Big) \dxdt                         \nonumber     \\
&    \quad \quad +  \int_0^{\tau} \int_{\Om}   \left(  1-\f{\vr }{  \tilde{\vr}  }     \right) p'(\tilde{\vr}) (\vv-\tilde{\vv}) \cdot \Grad \tilde{\vr} \dxdt
    +  \int_0^{\tau} \int_{\Om}          (\tilde{b}-b)        (\vv-\tilde{\vv}) \cdot \Grad \tilde{b} \dxdt         \nonumber     \\
    &   \quad \quad   -   \int_0^{\tau} \int_{\Om}    \Big(  p(\vr)-p(\tilde{\vr})  -p'(\tilde{\vr})(\vr-\tilde{\vr})    \Big) \Div \tilde{\vu}  \dxdt          \label{weak-strong-12}    \\
    & \quad \quad    -    \int_0^{\tau} \int_{\Om}    \Big(  \f{1}{2} b^2 -\f{1}{2} |\tilde{b}|^2 -\tilde{b}(b-\tilde{b} )                          \Big) \Div \tilde{\vu}   \dxdt  .          \nonumber
\end{align}

\subsection{Gronwall-type argument}\label{gronwall}

We are left to estimate the right-hand side of (\ref{weak-strong-12}) so as to invoke the Gronwall inequality. First of all, it is easily seen that
\begin{align}
    &     \left| \int_0^{\tau} \int_{\Om}    \Big(  p(\vr)-p(\tilde{\vr})  -p'(\tilde{\vr})(\vr-\tilde{\vr})    \Big) \Div \tilde{\vu}  \dxdt  \right|        \nonumber       \\
    &   \quad  \leq C         \int_0^{\tau} \int_{\Om}\Big(  H(\vr)-H(\tilde{\vr})-H'(\tilde{\vr})(\vr-\tilde{\vr}) \Big) \dxdt,        \nonumber   \\
    & \left|     \int_0^{\tau} \int_{\Om}    \Big(  \f{1}{2} b^2 -\f{1}{2} |\tilde{b}|^2 -\tilde{b}(b-\tilde{b} )                          \Big) \Div \tilde{\vu}   \dxdt           \right|  \label{weak-strong-13}   \\
    & \quad   \leq C      \int_0^{\tau} \int_{\Om} \f{1}{2}|b- \tilde{b}|^2      \dxdt .          \nonumber
\end{align}
Next, it is also clear that
\beq\label{weak-strong-14}
\left|  \int_0^{\tau} \int_{\Om}   \vr (\vv-\tilde{\vv} ) \cdot \Grad \tilde{\vu}  \cdot (\tilde{\vv}-\vv )  \dxdt \right| \leq
C \int_0^{\tau} \int_{\Om}    \f{1}{2}\vr |\vv-\tilde{\vv}|^2      \dxdt.
\eeq
One readily checks that the remaining terms are controlled by
\[
\int_0^{\tau} \int_{\Om}    |b-\tilde{b}| |\vv-\tilde{\vv}|  \dxdt,\,\, \int_0^{\tau} \int_{\Om}    |\vr-\tilde{\vr}| |\vv-\tilde{\vv}|  \dxdt.
\]
For instance,
\begin{align}
   \int_0^{\tau} \int_{\Om}    |b-\tilde{b}| |\vv-\tilde{\vv}|  \dxdt  &    \leq \ep \int_0^{\tau} \int_{\Om}   |\vv-\tilde{\vv}|^2     \dxdt +\f{C}{\ep}\int_0^{\tau} \int_{\Om} |b-\tilde{b}|^2 \dxdt   \nonumber  \\
    &   \leq   \ep   \int_0^{\tau}  \|  \vv-\tilde{\vv} \|_{W_{0}^{1,2}}^2 \dt +   \f{C}{\ep} \int_0^{\tau}\int_{\Om} |b-\tilde{b}|^2 \dxdt  \label{weak-strong-15}   \\
    &   \leq  C  \ep     \int_0^{\tau} \int_{\Om}          \mathbb{S}(\Grad (\vv-\tilde{\vv})):\Grad (\vv-\tilde{\vv}) \dxdt  +   \f{C}{\ep} \int_0^{\tau} \int_{\Om}|b-\tilde{b}|^2 \dxdt  ,       \nonumber
\end{align}
where the Korn-Poincar\'{e} inequality was used in the last step. Observe that the other integral $\int_0^{\tau} \int_{\Om}    |\vr-\tilde{\vr}| |\vv-\tilde{\vv}|  \dxdt$ can be estimated in the same spirit, up to a decomposition of the domain into its essential part and residual part. This is nowadays a routine matter and we may refer to Section 4 in Feireisl et al. \cite{FJA12} for the details.

Combining (\ref{weak-strong-13})-(\ref{weak-strong-15}) and fixing $\ep>0$ sufficiently small, we conclude from (\ref{weak-strong-12}) that
\[
\mathcal{E}(\vr,\vu,b\,|\, \tilde{\vr},\tilde{\vu},\tilde{b}  ) (\tau)  \leq   C    \int_0^{\tau}\mathcal{E}(\vr,\vu,b\,|\, \tilde{\vr},\tilde{\vu},\tilde{b}  )(t) \dt,
\]
yielding the desired result of Theorem \ref{TH2} via Gronwall inequality.

\section{Conclusion and remarks}\label{con-rem}
We consider a two-dimensional compressible viscous non-resistive MHD equations with general inflow-outflow boundary conditions. Such a general framework facilitates considerably the applications in astrophysics and thermonuclear reactions, among many others. Existence of global-in-time weak solutions is proved
for a fairly large class of physically admissible initial data. To our best knowledge, this is the first rigorous mathematical result concerning \emph{global-in-time} solvability of compressible viscous non-resistive MHD equations with \emph{general} boundary conditions. Based on the relative entropy method, we furthermore prove that a weak solution coincides with the strong solution if they share the same initial-boundary data. The weak-strong uniqueness property is also expected to be useful when studying the convergence of numerical solutions.

Several remarks concerning the main results are given below.
\begin{itemize}
\item{
In a recent work of Wen \cite{W21}, the author proved the existence of weak solutions to the compressible two-fluid model (\ref{in5})-(\ref{in6}) without any domination condition between the initial densities. As pointed out in Remark 1.4 of \cite{W21}, the method therein could be applied to our compressible viscous non-resistive MHD system (\ref{in2}) without any domination condition between the initial density and magnetic field. However, in our context of inflow-outflow boundary conditions, it seems impossible to remove the domination condition for \emph{any} $\gamma \geq \f{3}{2}$. This is due to the fact that the improved integrability estimates for the density hold only on compact subsets of the domain, cf. (\ref{pass-10}). Our main Theorem \ref{TH1} exploits the specific mathematical structure from the magnetic field and holds for \emph{any} $\gamma>1$.
}

\item{
To simplify the presentation, we have chosen the isentropic law. Nevertheless, in view of the methods from \cite{CBN19,Fei02}, it is possible to prove our main Theorem \ref{TH1} for a general non-monotone pressure-density law
obeying
\begin{align*}
    &   p(\vr)=\pi(\vr)  -\mathfrak{p} (\vr) ,\, \pi(\cdot) \in C([0,\infty)) \cap C^1((0,\infty)) ,\, \pi(0)=0,     \\
    &   \pi'(\vr) > \max\{  0, a_1\vr^{\gamma-1}-b \} ,\, \pi(\vr) \leq a_2 \vr^{\gamma}+a_3\,\, \text{  for any } \vr>0,       \\
    &   \mathfrak{p} (\cdot) \in C_c^2 ([0,\infty)) ,\, \mathfrak{p} (\cdot) \geq 0,\,  \mathfrak{p}(0)= \mathfrak{p}'(0)=0,
\end{align*}
for some $\gamma>1$ and positive constants $a_1,a_2,a_3$.
}


\item{ It is unknown whether the approximation scheme (\ref{appro-conti})-(\ref{appro-magn-1}) admits a unique global strong solution provided the initial-boundary data are suitably regular and the domain is smooth.
The difficulties are two-fold: the De Giorgi-Nash-Moser estimate for parabolic equation with general boundary conditions and a suitable control of velocity belonging to $L^q(0,T;L^r(\Om;\R^2)  )$ with $\f{2}{q}+ \f{2}{r}<1$.
In the context of homogeneous Dirichlet boundary conditions, this is achieved and we refer to Proposition 2.1 in \cite{LS19}.
}
\item{
Regularity of strong solutions in Theorem \ref{TH2} is certainly not optimal. It could be relaxed through a density argument; we shall not pursue this here.
}
\end{itemize}

\vspace{5mm}

\centerline{\bf Acknowledgements}
\vspace{2mm}
The research of Y. Li was supported by Natural Science Foundation of Anhui Province under grant number 2408085MA018, Natural Science Research Project in Universities of Anhui Province under grant numbers 2024AH050055 and 2023AH050105.
The work of Y.-S. Kwon was partially supported by the National Research Foundation of Korea (NRF2022R1F1A1073801).
The research of Y. Sun is supported by National Natural Science Foundation of China under grant number 12071211.
Last but not least, the authors sincerely thank the anonymous reviewers for many helpful suggestions.

\vspace{5mm}

\centerline{\bf Conflict of interest}
\vspace{2mm}
On behalf of all authors, the corresponding author states that there is no conflict of interest.

\vspace{5mm}

\centerline{\bf Data Availability Statement}
\vspace{2mm}
Data sharing not applicable to this article as no datasets were generated or
analyzed during the current study.


\end{document}